\documentclass[reqno, 11pt]{amsart}
\usepackage{amssymb,amsfonts, amsmath}
\usepackage{mathrsfs}
\usepackage{color}
\usepackage{dsfont}
\usepackage[normalem]{ulem} 
\newtheorem{theorem}{Theorem}[section]

\newtheorem{lemma}{Lemma}[section]

\numberwithin{equation}{section}

\renewcommand{\b}[1]{\mbox{\boldmath $#1$}}
\newcommand{\diag}{\mathrm{diag}}

\newcommand{\argmax}{\mathrm{argmax}}

\theoremstyle{definition}

\newtheorem{remark}{Remark}

\definecolor{gray}{rgb}{0.9,0.9,0.9}

\input cyracc.def

\begin{document}

\title[]{Moderate deviations for a diffusion type process in random environment}
\author{P. Chigansky}
\address{Universite du Maine,
Faculte des Sciences, Laboratoire de Statistique et Processus,
Avenue Olivier Messiaen, 72085 Le Mans CEDEX 9}
\email{Pavel.Chigansky@univ-lemans.fr}
\thanks{The research of P. Chigansky is supported by the Chateaubriand fellowship}

\author{R. Liptser}
\address{Department of Electrical Engineering Systems,
Tel Aviv University, 69978 Tel Aviv, Israel}
\email{liptser@eng.tau.ac.il}

\maketitle
\begin{abstract}
Let $\sigma(u)$, $u\in \mathbb{R}$ be an ergodic stationary Markov
chain, taking a finite number of values $a_1,\dots,a_m$, and
$b(u)=g(\sigma(u))$, where $g$ is a bounded and measurable function.

We consider the diffusion type process
$$
dX^\varepsilon_t = b(X^\varepsilon_t/\varepsilon)dt +
\varepsilon^\kappa\sigma\big(X^\varepsilon_t/\varepsilon\big)dB_t, \
t\le T
$$
subject to $X^\varepsilon_0=x_0$, where $\varepsilon$ is a small
positive parameter, $B_t$ is a Brownian motion, independent of
$\sigma$, and $\kappa> 0$ is a fixed constant. We show that for
$\kappa<1/6$, the family $\{X^\varepsilon_t\}_{\varepsilon\to 0}$
satisfies the Large Deviations Principle (LDP) of the
Freidlin-Wentzell type  with the constant drift $\mathbf{b}$ and the
diffusion  $\mathbf{a}$, given by
$$
\mathbf{b}=\sum\limits_{i=1}^m\dfrac{g(a_i)}{a^2_i}\pi_i\Big/
\sum\limits_{i=1}^m\dfrac{1}{a^2_i}\pi_i, \quad
\mathbf{a}=1\Big/\sum\limits_{i=1}^m\dfrac{1}{a^2_i}\pi_i,
$$
where $\{\pi_1,\ldots,\pi_m\}$ is the invariant distribution of the
chain $\sigma(u)$.
\end{abstract}

\maketitle
\section{Introduction}
\label{sec-1} M. Freidlin and R. Sowers, \cite{FrSow}, study the LDP
for the vector diffusion Markov processes, defined by the It\^o
equations with respect to Brownian motion $B_t$,
\begin{equation}\label{00}
X^\varepsilon_t=x_0+\int_0^tb(X^\varepsilon_s/\varepsilon)ds+\varepsilon^\kappa
\int_0^t\sigma(X^\varepsilon_s/\varepsilon)dB_s,
\end{equation}
where $\varepsilon>0$ is a small parameter and $b(x)$ and
$\sigma(x)$ are smooth periodic functions with period $1$. The
existence of three different LDP regimes, depending on the value of
$\kappa$, is shown in \cite{FrSow}:
$$
\kappa\in\Big(0,\frac{1}{2}\Big),\quad
\kappa=\frac{1}{2}\quad\text{and} \quad
\kappa\in\Big(\frac{1}{2},\infty\Big),
$$
where in the second and the third regimes the rate functions are not
of the classic LDP of Freidlin-Wentzell's type. The first regime is
characterized by the same rate function as for a diffusion process
$\widehat{X}^\varepsilon_t$ with constant drift and diffusion
parameters. Particularly, in the scalar case
\begin{equation*}
\widehat{X}^\varepsilon_t=x_0+\mathbf{b}t+\varepsilon^\kappa\sqrt{\mathbf{a}}B_t,
\end{equation*}
where $
\mathbf{b}=\int_0^1\frac{b(s)ds}{\sigma^2(s)}\big/\int_0^1\frac{ds}{\sigma^2(s)}
\quad\text{and}\quad
\mathbf{a}=1\big/\int_0^1\frac{1}{\sigma^2(s)}ds. $ In this context,
following the terminology of \cite{FW}, we shall refer LDP for
$\kappa\in(0,\frac 1 2)$ as MDP (Moderate Deviations Principle).

The aim of this note is to extend the MDP to a scalar diffusion in a
random environment, namely, when $b(u)=b(\omega,u)$ and
$\sigma(u)=\sigma(\omega,u)$ are random processes, independent of
the Brownian motion.

We assume that $\sigma(u)$ is a stationary and ergodic Markov chain
with a finite alphabet
$$
\mathcal{A}=\{a_1,\ldots,a_m\}, \ a_i\ne 0,
$$
having right continuous paths with left limits. We assume that
$$
b(u)=g(\sigma(u))
$$
for some bounded measurable function $g(x)$.

We prove the existence of MDP in the random environment for
\begin{equation}\label{6}
\kappa\in\Big(0, \frac{1}{6}\Big) \ \text{and} \
\mathbf{b}=\sum\limits_{i=1}^m\dfrac{g(a_i)}{a^2_i}\pi_i\Big/
\sum\limits_{i=1}^m\dfrac{1}{a^2_i}\pi_i, \quad
\mathbf{a}=1\Big/\sum\limits_{i=1}^m\dfrac{1}{a^2_i}\pi_i,
\end{equation}
where $\pi=(\pi_1,\ldots,\pi_m)$ is the invariant distribution of
$\sigma(u)$.

In \cite{FrSow} the LDP is derived using the Markov property of
$X^\varepsilon_t$ and asymptotic analysis, as $\varepsilon\to 0$, of
the log moment generating function (hereafter $\lambda(t)$ is a test
function):
\begin{equation}\label{limit}
\lim_{\varepsilon\to 0}\varepsilon^{2\kappa}\log
\mathsf{E}\Big(\exp\Big[\frac{1}
{\varepsilon^{2\kappa}}\int_0^T\lambda(t)dX^\varepsilon_t\Big]\Big)
\end{equation}
(see, e.g., Ch. 2.3 and Ch. 5.1 in Dembo and Zeitouni, \cite{DZ},
Dupuis and  Ellis, \cite{D+E}, Gartner and Freidlin,
\cite{GFe}).
When the environment is random, $X^\varepsilon_t$ is not a Markov
process anymore and calculation of the limit \eqref{limit} is quite
involved. Instead of \eqref{limit}, we apply the Puhalskii approach
based on the martingale exponential
$$
\mathfrak{z}_T(\lambda)=\exp\Big[\frac{1}
{\varepsilon^{2\kappa}}\int_0^T\lambda(t)dX^\varepsilon_t-
\log\mathscr{E}_T(X^\varepsilon)\Big]
$$
analysis, where
\begin{equation*}
\mathscr{E}_t(X^\varepsilon)=\exp\Big(\frac{1}{\varepsilon^{2\kappa}}
\int_0^t\Big[\lambda(s)b(X^\varepsilon_s/\varepsilon)-\frac{
\lambda^2(s)}{2}\sigma^2(X^\varepsilon_s/\varepsilon)\Big]ds\Big)
\end{equation*}
is the {\em cumulant} process compensating $ \exp\big(\frac{1}
{\varepsilon^{2\kappa}}\int_0^T\lambda(t)dX^\varepsilon_t\big) $ up
to a local martingale. Sufficient conditions for the MDP can be
formulated in terms of the cumulant process (see Corollary 4.3.8  in
\cite{puh2}) and in our setting even directly in terms of the drift
and diffusion parameters: by Corollary 6.7 in \cite{puh-ssr} the MDP
holds if for any $\eta>0$ one can choose $\mathbf{b}$,
$\mathbf{a}>0$ and  $\kappa>0$ such that
\begin{equation}\label{1f}
\begin{aligned}
& \lim_{\varepsilon\to
0}\varepsilon^{2\kappa}\log\mathsf{P}\Big(\sup_{t\le T}
\Big|\int_0^t[b(X^\varepsilon_s/\varepsilon)-\mathbf{b}]ds\Big|>\eta\Big)=-\infty
\\
& \lim_{\varepsilon\to
0}\varepsilon^{2\kappa}\log\mathsf{P}\Big(\sup_{t\le T}
\Big|\int_0^t[\sigma^2(X^\varepsilon_s/\varepsilon)-\mathbf{a}]ds\Big|>\eta\Big)=
-\infty;
\end{aligned}
\end{equation}
the appropriate choice is announced in \eqref{6}

The next section discusses the weak solutions of \eqref{00} in the
random environment. Our main result, Theorem \ref{theo-000}, is
formulated in Section \ref{sec-3} and is proved in Section
\ref{sec-5} which is preceded by auxiliary results in Section
\ref{sec-4}. The MDP gap between
$$
\kappa\in \Big(0,\frac{1}{2}\Big) \quad\text{and}\quad \kappa\in
\Big(0,\frac{1}{6}\Big)
$$
for oscillating and random environments  is discussed in Section
\ref{sec-6}.

For reader's fast reference,  the essential details of Puhalskii's
method, adapted to our setting, are outlined in Appendix.

\section{ Diffusion in random environment}
\label{sec-2}

Hereafter, we will deal with a weak solution of the scalar equation
explicitly constructed by time scaling and  change of probability
measure (for other approaches see \cite{Brox}, \cite{Schu}).

Let $\sigma=\big(\sigma(u)\big)_{u\in \mathbb{R}}$ be the Markov
chain, defined in the previous section, and $\beta=(\beta_t)_{t\ge
0}$ be a Brownian motion independent of $\sigma$. Assume that the
pair $(\sigma,\beta)$ is defined on a stochastic basis
$(\varOmega,\mathcal{F},\mathbf{F}=(\mathcal{F}_t)_{t\ge
0},\mathsf{Q})$ with the general conditions, $\sigma$ is
$\mathcal{F}_0$-measurable and $\beta$ is independent of
$\mathcal{F}_0$.

For $t>0$, introduce the stopping time $
\tau_t=\inf\big\{r:\int_0^{r}\frac{1}{\varepsilon^{2\kappa}
\sigma^2((\beta_s+x_0)/\varepsilon)}ds\ge t\big\}. $ Since
$\sigma^2(u)>0$, we have $
\int_0^{\tau_t}\frac{1}{\varepsilon^{2\kappa}
\sigma^2((\beta_s+x_0)/\varepsilon)}ds\equiv t $ and, in turn,
$$
\tau_t=\int_0^t\varepsilon^{2\kappa}\sigma^2\big((\beta_{\tau_s}+x_0)/
\varepsilon\big)ds.
$$
Introduce the filtration $\mathbf{G}=(\mathcal{G}_t)_{t\ge 0}$ with
$\mathcal{G}_t:=\mathcal{F}_{\tau_t}$. Obviously,
$(\beta_{\tau_t},\mathcal{G}_t)$ is a continuous martingale with the
quadratic variation process $\tau_t$. Then by the Levy-Doob theorem
the process
$$
B_t= \int_0^t\frac{1}{\varepsilon^{\kappa}\sigma\big((\beta_{\tau_s}
+x_0)/\varepsilon\big)}d\beta_{\tau_t}
$$
is Brownian motion. Since $B=(B_t)_{t\ge 0}$ is independent of
$\mathcal{F}_0$ and $\mathcal{F}_0=\mathcal{G}_0$, $B=(B_t)_{t\ge
0}$ is independent of $\mathcal{G}_0$. On the other hand, by the
same reason $\sigma$ is $\mathcal{G}_0$-measurable. Hence, $(B_t)$
and $\sigma(u)$ are independent random processes.

Denote $Y_t:= x_0 + \beta_{\tau_t}$. Then, the definition of $B_t$
implies the following representation for $Y_t$:
\begin{equation}\label{YY}
Y_t=x_0+\int_0^t\varepsilon^{\kappa}\sigma(Y_s/\varepsilon)dB_s.
\end{equation}
Consequently, at least one weak solution of  \eqref{00} with zero
drift exists.
 A weak solution of \eqref{00} with the required drift
can be constructed with the help of Girsanov's theorem. With $Y_t$,
defined  in \eqref{YY}, set
$$
\Upsilon_T=\exp\left(\int_0^T\frac{b(Y_s/\varepsilon)}
{\varepsilon^\kappa\sigma(Y_s/\varepsilon)}dB_s
-\frac{1}{2}\int_0^T\frac{b^2(Y_s/\varepsilon)}{\varepsilon^{2\kappa}
\sigma^2(Y_s/\varepsilon)}ds \right).
$$
Since
$\dfrac{b(Y_s/\varepsilon)}{\varepsilon^\kappa\sigma(Y_s/\varepsilon)}$
is bounded and $T<\infty$, we have $
\int_\varOmega\Upsilon_Td\mathsf{Q}=1. $ We define a probability
measure $\mathsf{P}$ with $d\mathsf{P}:=\mathfrak{z}_Td \mathsf{Q}$.
Then, by the Girsanov theorem,
$$
\widehat{B}_t=B_t-\int_0^t
\frac{b(Y_s/\varepsilon)}{\varepsilon^\kappa\sigma(Y_s/\varepsilon)}ds
$$
is the Brownian motion with respect to $\mathbf{G}$ under
$\mathsf{P}$. In other words, the process $Y_t$ defined on the new
stochastic basis
$\big(\varOmega,\mathcal{F},\mathbf{G}=(\mathcal{G}_t)_{t\ge
0},\mathsf{P}\big)$ admits the following representation
\begin{align*}
Y_t=x_0+\int_0^tb(Y_s/\varepsilon)ds+\int_0^t
\varepsilon^\kappa\sigma(Y_s/\varepsilon)d\widehat{B}_s,
\end{align*}
that is, $Y_t$ is a weak solution of \eqref{00}. Since
$(\widehat{B}_t)_{t\le T}$ is $\mathsf{P}$-independent of
$\mathcal{G}_0$ and $(\sigma(u))_\mathbb{R}$ is
$\mathcal{G}_0$-measurable, the Brownian motion $\widehat{B}_t$ and
$\sigma(u)$ are $\mathsf{P}$-independent random processes.

\section{The main result}
\label{sec-3}

Let $X^\varepsilon=(X^\varepsilon_t)_{t\le T}$ be a weak solution of
\eqref{00}. Recall that $X^\varepsilon$ satisfies LDP (in our case
MDP) with the rate $\varepsilon^{2\kappa}$ and the good rate
function $J(u)$ in the space of continuous functions
$\mathbb{C}_{[0,T]}$ endowed with the uniform metric if for any
closed set $F$ and open set $G$,
\begin{align*}
& \varlimsup_{\varepsilon\to
0}\varepsilon^{2\kappa}\log\mathsf{P}(X^\varepsilon \in F)\le
-\inf_{u \in F}J(u)
\\
& \varliminf_{\varepsilon\to
0}\varepsilon^{2\kappa}\log\mathsf{P}(X^\varepsilon \in G)\ge
-\inf_{u \in G}J(u).
\end{align*}

\begin{theorem}\label{theo-000}
For $\kappa<1/6$, the family $\{X^\varepsilon\}_{\varepsilon\to 0}$
satisfies the MDP   with the rate $\varepsilon^{2\kappa}$ and the
rate function
$$
J(u)=
  \begin{cases}
    \frac{1}{2\mathbf{a}}\int_0^T[\dot{u}_t-\mathbf{b}]^2dt , &
\begin{subarray}
\phantom{}u_0=x_0\\ du_t=\dot{u}_tdt
\\\int_0^T\dot{u}^2_tdt<\infty
\end{subarray}
    \\
  \infty , & \text{otherwise}
  \end{cases}
$$
with
\begin{equation*}
\mathbf{b}=\sum\limits_{i=1}^m\frac{g(a_i)\pi_i}{a^2_i}\Big/
\sum\limits_{i=1}^m\frac{\pi_i}{a^2_i} \quad\text{and}\quad
\mathbf{a}=1\Big/ \sum\limits_{i=1}^m\frac{\pi_i}{a^2_i},
\end{equation*}
where $\{\pi_1,...,\pi_m\}$ is the invariant distribution of
$\sigma$.
\end{theorem}
The proof of this theorem requires some auxiliary results gathered
in the next section.
\section{Auxiliary results}
\label{sec-4}

Henceforth,

- $A^*$ is transposed of a matrix $A$;

- for any $x\in\mathbb{R}^d$, $\diag(x)$ is the diagonal matrix with
(diagonal) entries $x_i$'s;

- $l$ is a generic positive constant whose meaning may change from
line to line;

- $\inf\{\varnothing\}=\infty$.

\subsection{\bf The Poisson decomposition}
\label{sec-4.1}

Let $\mathbf{F}^\sigma=(\mathcal{F}^\sigma_t)_{t\in\mathbb{R}}$ be
the filtration generated by $\sigma$:
$\mathcal{F}^\sigma_t=\{\sigma(u), -\infty<u\le t\}$. Since $\sigma$
is an ergodic chain, its transition intensities matrix $\varLambda$
has simple zero eigenvalue. Therefore, for any bounded measurable
function $\Psi(x)$ with $ \mathsf{E} \Psi\big(\sigma(0)\big)=0$
there exists $\gamma>0$ such that
$|\mathsf{E}\big(\Psi(\sigma(t))|\mathcal{F}^\sigma_0 \big)|\le
le^{-\gamma t}$ a.s. for any $t>0$. Hence, $
\int_0^\infty\mathsf{E}|\mathsf{E}\big(\Psi(\sigma(t)\big|\mathcal{F}^\sigma_0\big)
\big|dt<\infty. $ Then (see e.g. Ch 9, \S 2 in \cite{LSMar}) the
process $\int_0^t\Psi(\sigma(s))ds$ obeys the Poisson decomposition
\begin{equation}\label{VM}
\int_0^t\Psi(\sigma(s))ds=V_t-V_0-M_t,
\end{equation}
where $V_t$ is $\mathbf{F}^\sigma$-adapted process and $M_t$ is
$\mathbf{F}^\sigma$-martingale with right continuous pathes having
left limits. In the case under consideration, $M_t$ is a square
integrable martingale (see Lemma \ref{lem-2.1}) with the quadratic
variation process $\langle M\rangle_t$.

\begin{lemma}\label{lem-2.1}
\mbox{}

{\rm 1)} $|V_t|\le l$ for any $t\ge 0${\rm ;}

{\rm 2)} $M_t$ is a purely discontinuous square integrable
martingale with bounded jumps{\rm ;}

{\rm 3)} $d\langle M\rangle_t=m(t)dt$, $m(t)\le l$.
\end{lemma}

\begin{proof}
Denote by
$$
I(t)=
\begin{pmatrix}
  I_{\{\sigma(t)=a_1\}} \\
  \vdots \\
  I_{\{\sigma(t)=a_m\}} \\
\end{pmatrix}
\quad \text{and}\quad f=
\begin{pmatrix}
  \Psi(a_1) \\
  \vdots \\
  \Psi(a_m) \\
\end{pmatrix}.
$$
The obvious equality $\Psi (\sigma(t))=f^*I(t)$ implies
$$
0=\mathsf{E}\Psi\big(\sigma(t)\big)=f^*\mathsf{E}I(t)=f^*\pi.
$$
This property of $f$ and the aforementioned spectral gap of the
matrix $\varLambda$ guarantees solvability of the Poisson equation
\begin{equation}\label{lae}
\varLambda g=f
\end{equation}
whose solution is unique in the class $g^*\pi=0$. Only this solution
will be considered in the sequel.

By Lemma 9.2, Ch.9, \S 9.1 in \cite{LSI},
\begin{equation}\label{Nt}
N_t:=I(t)-I(0)-\int_0^t\Lambda^* I(s)ds
\end{equation}
is a purely discontinuous martingale, with respect to
$\mathbf{F}^\sigma$, with bounded jumps. We show now that
$V_t=g^*I(t)$ and $M_t=g^*N_t$. Multiplying from the left both sides
of \eqref{Nt} by $g^*$ and taking into account the definition of
$V_t$ and $M_t$ we find that $
M_t=V_t-V_0-\int_0^tg^*\varLambda^*I(s)ds. $ Further, by
\eqref{lae}, $ g^*\varLambda^*I(s)=f^*I(s)=\Psi(\sigma(u)). $ In
other words, \eqref{VM} holds true with $V_t$ and $M_t$ chosen
above. Therefore, statements 1) and 2) are obvious. The statement 3)
is proved as follows: by the It\^o formula we find that
\begin{align*}
I(t)I^*(t)&=I(0)I^*(0)+\int_0^t[I(s)I^*(s)\Lambda+\Lambda^*
I(s)I^*(s)]ds
\\
&\quad +\int_0^t[I(s-)dN^*_s+dN_sI^*(s-)]+[N,N]_t
\\
&=I(0)I^*(0)+\int_0^t[I(s)I^*(s)\Lambda+\Lambda^*
I(s)I^*(s)]ds+\langle N\rangle_t
\\
 &\quad +\text{martingale},
\end{align*}
where $\langle N\rangle_t$ is the quadratic variation process of
$N_t$, and, owing to the identity $I(t)I^*(t)=\diag(I(t))$, also
that
\begin{align*}
I(t)I^*(t)=I(0)I^*(0)+\int_0^t\diag\big(\Lambda^* I(s)\big)ds
+\text{martingale}.
\end{align*}
Both representations  for $I(t)I^*(t)$ imply that the  predictable
process with paths in the Skorokhod space  of locally bounded
variation
$$
\int_0^t\big([I(s)I^*(s)\Lambda+\Lambda^* I(s)I^*(s)]
-\diag[\Lambda^* I(s)]\big)ds+\langle N\rangle_t
$$
is a martingale starting from zero. Hence, by Theorem 1 in Ch. 2, \S
2, \cite{LSMar}, this martingale is indistinguishable from zero or,
equivalently,
$$
\langle N\rangle_t=\int_0^t\big(\diag[\Lambda^*
I(s)]-[I(s)I^*(s)\Lambda+\Lambda^* I(s)I^*(s)] \big)ds.
$$
Therefore, $ d\langle M\rangle_t\equiv g^*d\langle
N\rangle_tg=g^*\big(\diag[\Lambda^*
I(t)]-[I(t)I^*(t)\Lambda+\Lambda^* I(t)I^*(t)] \big)gdt. $
\end{proof}

\subsection{\bf Exponential estimate for martingales
with bounded jumps}

For a continuous martingale $M=(M_t)_{t\ge 0}$ with $M_0=0$ and the
quadratic variation process $\langle M\rangle_t$ the following
exponential estimate is well known (see e.g.  Lemma 1 in
\cite{LipSpok}): for any $q,r>0$,
\begin{align}\label{ccc}
 \mathsf{P}\Big(\sup_{t\le T}|M_t|\ge
 r , \ \langle M\rangle_T\le q\Big) \le
2\exp\Big(-\frac{ r ^2}{2q}\Big).
\end{align}
A similar inequality holds for discontinuous martingales.
\begin{lemma}\label{lem-L}
Let $M=(M_t)_{t\ge 0}$ be a purely discontinuous martingale with
$M_0=0$ and paths in the Skorokhod space $\mathbb{D}$ with  bounded
jumps $|\triangle M_t|\le K$ and the quadratic variation process
$\langle M\rangle_t$. Then, for any  $q, r>0$
\begin{equation}\label{D}
\mathsf{P}\Big(\sup_{t\le T}|M_t|\ge  r , \ \langle M\rangle_T\le
q\Big) \le 2\exp\Big(-\frac{ r ^2}{2(K r +q)}\Big).
\end{equation}
\end{lemma}
\begin{remark}
\eqref{ccc} is a particular case of \eqref{D} for $K=0$.
\end{remark}
\begin{proof}
Denote  by $\mu=\mu(dt,dz)$ the integer-valued measure,  associated
with the jump process $\triangle M_t$, and by $\nu=\nu(dt,dz)$ the
compensator of $\mu$ (for more details, see e.g. \cite{LSMar} or
\cite{JSn}); clearly $\nu(\mathbb{R}_+\times\{|z|>K\})=0$ since
$|\triangle M_t|\le K$.

Then $ M_t=\int_0^t\int_{|z|\le K}z[\mu(ds,dz)-\nu(ds,dz)] $ and $
\langle M\rangle_t= \int_0^t\int_{|z|\le K}z^2\nu(ds,dz). $ Let
$\mathscr{L}_t(\lambda)$ be the cumulant process associated with
$M_t$, i.e. for any $\lambda\in\mathbb{R}$ the random process
$$
\mathfrak{z}_t(\lambda)=e^{\lambda M_t-\mathscr{L}_t(\lambda)}
$$
is a local martingale. It can be easily checked with the help of
It\^o's formula that
$$
\mathscr{L}_t(\lambda)=G_t(\lambda)+\sum_{s\le
t}\big[\log(1+\triangle G_s(\lambda)) -\triangle G_s(\lambda)\big],
$$
where
\begin{align*}
G_t(\lambda)&=\int_0^t\int_{|z|\le K} \big(e^{\lambda z}-1-\lambda
z\big)\nu(ds,dz)
\\\triangle G_t(\lambda)&=\int_{|z|\le K}\big(e^{\lambda z}-1-\lambda
z\big) \nu(\{t\},dz).
\end{align*}
The positive local martingale $\mathfrak{z}_t(\lambda)$ is also a
supermartingale. Hence $ \mathsf{E}\mathfrak{z}_\tau(\lambda)\le 1 $
for any stopping time $\tau$. Since $\triangle G_t(\lambda)\ge 0$
and  $\log(1+x)-x\le 0$ for $x\ge 0$, we have
$\mathscr{L}_t(\lambda)\le G_t(\lambda)$. Consequently,
 $\mathfrak{z}_\tau\ge e^{\lambda
M_\tau-G_\tau(\lambda)}$ for any stopping time $\tau$, that is,
\begin{equation}\label{zzz<1}
\mathsf{E}e^{\lambda M_\tau-G_\tau(\lambda)}\le 1.
\end{equation}
For $|z|\le K$ and $|\lambda|<1/K$, we have
\begin{align*}
e^{\lambda z}-\lambda z-1&=\sum_{j=2}^\infty\frac{(\lambda z)^j}{j!}
\le \frac{\lambda^2z^2}{2}\sum_{j=0}^\infty|\lambda z|^j
=\frac{\lambda^2z^2}{2}\frac{1}{1-|\lambda z|}\le
\frac{\lambda^2z^2}{2}\frac{1} {1-|\lambda K|}.
\end{align*}
Hence
$$
G_\tau (\lambda)\le \int_0^\tau\int_{|z|\le
K}\frac{\lambda^2z^2}{2}\frac{1} {1-|\lambda K|}\nu(ds,dz)
=\frac{1}{1-|\lambda|K}\frac{\lambda^2}{2}\langle M\rangle_\tau.
$$

Now, due to \eqref{zzz<1}, for any measurable set $A$ we obtain that
$$
1\ge \mathsf{E}I_{\{A\}}\exp\Bigg(\lambda
M_\tau-\frac{1}{1-|\lambda|K}\frac{\lambda^2}{2} \langle
M\rangle_\tau\Bigg).
$$
The choice of $ \tau=\inf\{t\le T:M_t\ge r\} $  and $A=\{\tau\le
T\}\cap\{\langle M\rangle_T\le q\}$ for any $|\lambda|\le 1/K$
implies
\begin{align*}
1\ge \mathsf{E}I_{\{A\}} \exp\Bigg(\lambda
M_\tau-\frac{1}{1-|\lambda|K}\frac{\lambda^2}{2} \langle
M\rangle_\tau\Bigg)\ge \mathsf{E}I_{\{A\}} \exp\Bigg(\lambda  r
-\frac{1}{1-|\lambda|K}\frac{\lambda^2}{2} q\Bigg).
\end{align*}
So, taking into account that $\{\sup_{t\le T}M_t\ge  r \}=\{\tau\le
T\}$, we find that
$$
\mathsf{P}\Big(\sup_{t\le T}M_t\ge r ,\langle M\rangle_T\le q\Big)
\le \exp\Bigg(-\frac{\lambda  r (1-\lambda K)-
\frac{\lambda^2}{2}q}{1-\lambda K}\Bigg).
$$
Finally, the choice of $ \lambda'=\mathop{\argmax}
\limits_{\lambda\in(0,\frac{1}{K})}\big[\lambda  r (1-\lambda K) -
\frac{\lambda^2}{2}q\big] $ provides
$$
\mathsf{P}\big(\sup_{t\le T}M_t\ge  r , \ \langle M\rangle_T\le
q\big) \le \exp\bigg(-\frac{ r ^2}{2(K r +q)}\bigg).
$$
The same inequality holds for $\sup_{t\le T}(-M_t)$.

Now, \eqref{D} follows from $ \mathsf{P}(A\cup B)\le
2\big[\mathsf{P}(A)\vee\mathsf{P}(B)\big], $ for any measurable sets
$A$ and $B$.
\end{proof}

\section{The proof of Theorem \ref{theo-000}} \label{sec-5}

Recall that \eqref{1f} implies the required MDP. We begin with the
proof of the first part in \eqref{1f}. Introduce the stationary
process
\begin{equation*}
\theta(t)=\frac{b(t)-\mathbf{b}}{\sigma^2(t)},
\end{equation*}
where   $\mathbf{b}$ is a fixed constant such that
$\mathsf{E}\theta(t)=0$; in other words, $
\mathbf{b}=\sum\limits_{i=1}^m\frac{g(a_i)\pi_i}{a^2_i}\big/
\sum\limits_{i=1}^m\frac{\pi_i}{a^2_i}. $

Define
\begin{equation*}
H(x)=\int_0^x\int_0^v\theta (s)dsdv.
\end{equation*}
The random function $H(x)$ is continuously differentiable and has
bounded Sobolev's second derivative. Hence by Krylov's version of
the It\^o formula \cite{KI}
\begin{align*}
H(X^\varepsilon_t/\varepsilon)&=H(x_0/\varepsilon)+
\frac{1}{\varepsilon}\int_0^t\int_0^{X^\varepsilon_v/\varepsilon}
\theta(s)dsdX^\varepsilon_v +\frac{1}{\varepsilon^{2(1-\kappa)}}
\int_0^t[b(X^\varepsilon_s/\varepsilon)- \mathbf{b}]ds
\end{align*}
or, equivalently,
\begin{multline}\label{5.3a}
\int_0^t[b(X^\varepsilon_s/\varepsilon)-\mathbf{b}]ds =
\varepsilon^{2(1-\kappa)}\int_0^{X^\varepsilon_t/\varepsilon}\int_0^v
\theta(s)dsdv -
\varepsilon^{2(1-\kappa)}\int_0^{x_0/\varepsilon}\int_0^v
\theta(s)dsdv
-\\
\varepsilon^{1-2\kappa}\int_0^t\int_0^{X^\varepsilon_v/\varepsilon}
\theta(s) dsb(X^\varepsilon_v/\varepsilon) dv
-\varepsilon^{1-\kappa}\int_0^t\int_0^{X^\varepsilon_v/\varepsilon}
\theta(s)ds\sigma(X^\varepsilon_v/\varepsilon)dB_v.
\end{multline}
Let $Z^\varepsilon_t$ denote any of the terms in the right hand side
of \eqref{5.3a}. Obviously, the first part \eqref{1f} holds true if
\begin{equation}\label{11f}
 \lim_{\varepsilon\to 0}\varepsilon^{2\kappa}\log\mathsf{P}\Big(\sup_{t\le T}
\big|Z^\varepsilon_t\big|>\eta\Big)=-\infty.
\end{equation}
In order to simplify the proof of \eqref{11f}, let us show that the
set $\Upsilon^\varepsilon_C=\big \{\sup_{t\le
T}|X^\varepsilon_t|>C\big\}$ is exponentially negligible in the
sense
\begin{equation}\label{Ups}
\lim_{C\to\infty}\varlimsup_{\varepsilon\to
0}\varepsilon^{2\kappa}\log\mathsf{P}\big(
\Upsilon^\varepsilon_C\big)=-\infty.
\end{equation}
Denote by
$A^\varepsilon_t=x_0+\int_0^tb(X^\varepsilon_s/\varepsilon)ds$ and
$M^\varepsilon_t=\varepsilon^\kappa\int_0^t
\sigma(X^\varepsilon_s/\varepsilon)dB_s$. Since
$X^\varepsilon_t=A^\varepsilon_t+ M^\varepsilon_t$ and $\sup_{t\le
T} |A^\varepsilon_t|$ is bounded by a constant independent of
$\varepsilon$, the proof of \eqref{Ups} reduces to the proof of
$$
\lim_{C\to\infty}\varlimsup_{\varepsilon\to
0}\varepsilon^{2\kappa}\log\mathsf{P} \big(\sup_{t\le
T}\big|M^\varepsilon_t\big|>C\big)=-\infty.
$$
Since $ \mathsf{P}\big(\langle M^\varepsilon\rangle_T\le
lT\varepsilon^{2\kappa}\big)=1, $ due to \eqref{ccc}, we have $
\mathsf{P}\big(\sup_{t\le T}\big|M^\varepsilon_t\big|>C\big) \le
2\exp\big(-\frac{C^2}{2\varepsilon^{2\kappa}lT}\big) $ and in turn
\eqref{Ups}.

In view of \eqref{Ups}, instead of \eqref{11f} it suffices to prove
that for any $C>0$
\begin{equation}\label{only}
\varlimsup_{\varepsilon\to
\infty}\varepsilon^{2\kappa}\log\mathsf{P}\Big(\sup_{t\le
T}\big|Z^\varepsilon_t \big|\ge\eta,
\varOmega\setminus\Upsilon^\varepsilon_C\Big)=-\infty.
\end{equation}

Let
$Z^\varepsilon_t:=\varepsilon^{2(1-\kappa)}\int_0^{X^\varepsilon_t/\varepsilon}
\int_0^v\theta(s)dsdv. $ Then, $\sup_{t\le T}|Z^\varepsilon_t|$ is
bounded from above on the set
$\{\varOmega\setminus\Upsilon^\varepsilon_C\}$ by $ C\sup_{|v|\le
C}\big|\varepsilon^{1-2\kappa}\int_0^{v/\varepsilon}
\theta(s)ds\big|. $ Consequently,
\begin{align*}
\mathsf{P}\Big(\sup_{t\le T}\big|Z^\varepsilon_t \big|\ge\eta,
\varOmega\setminus\Upsilon^\varepsilon_C\Big) &\le
\mathsf{P}\Big(\sup_{|v|\le
C}\big|\varepsilon^{1-2\kappa}\int_0^{v/\varepsilon}
\theta(s)ds\big|\ge \eta,
\varOmega\setminus\Upsilon^\varepsilon_C\Big)
\\
&\le \mathsf{P}\Big(\sup_{|v|\le
C}\big|\varepsilon^{1-2\kappa}\int_0^{v/\varepsilon}
\theta(s)ds\big|\ge \eta\Big)
\end{align*}
and \eqref{only} holds provided that
\begin{equation}\label{ux0}
\lim_{\varepsilon\to
0}\varepsilon^{2\kappa}\log\mathsf{P}\Big(\sup_{|v|\le C}
\Big|\varepsilon^{1-2\kappa}\int_0^{v/\varepsilon}
\theta(s)ds\Big|\ge \eta\Big)=-\infty.
\end{equation}

For $
Z^\varepsilon_t:=\varepsilon^{1-2\kappa}\int_0^t\int_0^{X^\varepsilon_v/\varepsilon}
\theta(s) dsb(X^\varepsilon_v/\varepsilon) dv$ and
$Z^\varepsilon_t:=\varepsilon^{2(1-\kappa)}\int_0^{x_0/\varepsilon}\int_0^{v}
\theta(s)dsdv$,
 \eqref{ux0}
implies \eqref{only}.

\smallskip
Recall that $\theta(s)$ is a  strictly stationary process and,
therefore, the distributions of
$$
\sup_{0\le v\le C}\Big|\int_0^{v/\varepsilon}\theta(s)ds\Big|
\quad\text{and}\quad \sup_{0\le v\le
C}\Big|\int^0_{-v/\varepsilon}\theta(s)ds\Big|
$$
coincide. Hence,
\begin{align*}
& \mathsf{P}\Big(\sup_{|v|\le C}
\Big|\varepsilon^{1-2\kappa}\int_0^{v/\varepsilon}
\theta(s)ds\Big|\ge \eta\Big)
\\
&\le 2\Bigg[\mathsf{P}\Big(\sup_{0\le v\le C}
\Big|\varepsilon^{1-2\kappa}\int_0^{v/\varepsilon}
\theta(s)ds\Big|\ge \frac{\eta}{2}\Big)\bigvee
\mathsf{P}\Big(\sup_{0\le v\le C}
\Big|\varepsilon^{1-2\kappa}\int^0_{-v/\varepsilon}
\theta(s)ds\Big|\ge \frac{\eta}{2}\Big)\Bigg]
\\
&=2\mathsf{P}\Big(\sup_{0\le v\le C}
\Big|\varepsilon^{1-2\kappa}\int_0^{v/\varepsilon}
\theta(s)ds\Big|\ge \frac{\eta}{2}\Big).
\end{align*}

Thus, instead of \eqref{ux0}, we shall prove
\begin{equation}\label{uxx}
\lim_{\varepsilon\to
0}\varepsilon^{2\kappa}\log\mathsf{P}\Big(\sup_{0<v\le C}
\Big|\varepsilon^{1-2\kappa}\int_0^{v/\varepsilon}
\theta(s)ds\Big|\ge \frac{\eta}{2}\Big)=-\infty.
\end{equation}
Since $\theta(s)=\Psi(\sigma(s))$ and $\mathsf{E}\theta(s)=0$, by
Lemma \ref{lem-2.1} we have
$$
\varepsilon^{1-2\kappa}\int_0^{v/\varepsilon}\theta(s)ds
=\varepsilon^{1-2\kappa}[V_{v/\varepsilon}-V_0]+\varepsilon^{1-2\kappa}
M_{v/\varepsilon},
$$
where $V_t$ is a bounded process and $M_t$ is a purely discontinuous
martingale with bounded jumps $|\triangle M_s|\le l$ and $\langle
M\rangle_v\le lv $. Therefore, we shall prove only that for any
$\eta>0$,
\begin{equation}\label{oj}
\lim_{\varepsilon\to
0}\varepsilon^{2\kappa}\log\mathsf{P}\Big(\sup_{0<v\le C}
\big|\varepsilon^{1-2\kappa}M_{v/\varepsilon}\big|\ge
\eta\Big)=-\infty.
\end{equation}

By  \eqref{D} and $ \mathsf{P}\big(\langle
\varepsilon^{1-2\kappa}M\rangle_{C/\varepsilon}\le
lC\varepsilon^{1-4\kappa}\big)=1, $ we have the following upper
bound for $\kappa<\frac{1}{4}:$
\begin{equation*}
\mathsf{P}\Big(\sup_{0\le v\le
C}|\varepsilon^{1-2\kappa}M_{v/\varepsilon}|\ge \eta\Big)
 \le 2\exp\Bigg(-\frac{\eta^2}{2(l\eta\varepsilon^{1-2\kappa}+lC\varepsilon^{1-4\kappa})}
 \Bigg).
\end{equation*}
Consequently for $\kappa<\frac{1}{6}$,
\begin{equation*}
\varepsilon^{2\kappa}\log \mathsf{P}\Big(\sup_{0\le v\le
C}|\varepsilon^{1-2\kappa}M_{v/\varepsilon}|\ge \eta\Big) \le
\varepsilon^{2\kappa}\log 2
-\frac{\eta^2}{\varepsilon^{1-6\kappa}2l(\eta\varepsilon^2+C)}
\end{equation*}
and \eqref{oj} follows.

So, it is left to verify \eqref{11f} for
$Z^\varepsilon_t:=\varepsilon^{1-\kappa}\int_0^t\int_0^{X^\varepsilon_v/\varepsilon}
\theta(s)ds\sigma(X^\varepsilon_v/\varepsilon)dB_v. $ Since
$Z^\varepsilon_t$ is a continuous martingale with
$
\langle Z^\varepsilon\rangle_t=\varepsilon^{2(1-\kappa)}\int_0^t
\big(\int_0^{X^\varepsilon_v/\varepsilon} \theta(s)ds\big)^2\sigma^2
(X^\varepsilon_v/\varepsilon)dv
$
and $\sigma^2\le l$ we have
$$
\langle Z^\varepsilon\rangle_T\le Tl
\sup_{|v|\le C}\Big(\varepsilon^{1-\kappa}\int_0^{v/\varepsilon} \theta(s)ds\Big)^2,
\quad
\text{on the set $\varOmega\setminus\Upsilon^\varepsilon_C$}.
$$
on the set $\varOmega\setminus\Upsilon^\varepsilon_C$,
For $r>0$, write
\begin{align}\label{1/6}
&\mathsf{P}\Big(\sup_{t\le T}|Z^\varepsilon_t|>\eta,
\varOmega\setminus  \Upsilon^\varepsilon_C\Big)\le
\mathsf{P}\Big(\sup_{t\le T}|Z^\varepsilon_t|>\eta, \ \langle
Z^\varepsilon\rangle_T\le Tl\sup_{|v|\le C}
\Big|\varepsilon^{1-\kappa}\int_0^{v/\varepsilon}
\theta(s)ds\Big|^2\Big)
\nonumber
\\
&\le 2\Big\{\mathsf{P}\Big(\sup_{t\le T}|Z^\varepsilon_t|>\eta, \
\langle Z^\varepsilon\rangle_T\le 2lr^2\varepsilon^{2\kappa}\Big)
\bigvee 2\mathsf{P}\Big(\sup_{0\le v\le C}
\Big|\varepsilon^{1-\kappa}\int_0^{v/\varepsilon}
\theta(s)ds\Big|>r\varepsilon^\kappa\Big)\Big\}
\nonumber
\\
&\le
 4\exp\Big(\frac{-\eta^2}{4l r^2 \varepsilon^{2\kappa}}\Big)
\bigvee 4 \mathsf{P}\Big(\sup_{0\le v\le C}
\Big|\varepsilon^{1-2\kappa}\int_0^{v/\varepsilon}
\theta(s)ds\Big|>r\Big).
\end{align}
Since  $r$ in  \eqref{1/6} is an arbitrary positive parameter, we
have
\begin{multline*}
\varlimsup_{\varepsilon\to 0}\varepsilon^{2\kappa}\log
\mathsf{P}\Big(\sup_{t\le T}|Z^\varepsilon_t|>\eta,
\varOmega\setminus  \Upsilon^\varepsilon_C\Big)
\\
\le \varepsilon^{2\kappa}\log 4-\frac{\eta^2}{4lr^2}\bigvee\varlimsup_{\varepsilon\to
0}\varepsilon^{2\kappa} \log\mathsf{P}\Big(\sup_{0<v\le C}
\Big|\varepsilon^{1-2\kappa}\int_0^{v/\varepsilon}
\theta(s)ds\Big|>r\Big) =-\frac{\eta^2}{4lr^2}\xrightarrow[r\to
0]{}-\infty,
\end{multline*}
where the equality is due to \eqref{uxx} proved above.

\medskip
The proof for the second part in \eqref{1f} is similar: we introduce
the stationary process
\begin{equation*}
\theta(t)=1-\frac{\mathbf{a}}{\sigma^2(t)},
\end{equation*}
where   $\mathbf{a}$ is a fixed constant such that
$\mathsf{E}\theta(t)=0$, i.e. $
\mathbf{a}=1\big/\sum\limits_{i=1}^m\frac{\pi_i}{a^2_i}, $ and set $
H(x)=\int_0^x\int_0^v\theta (s)dsdv. $ By Krylov-It\^o's formula
\cite{KI},
\begin{align*}
H(X^\varepsilon_t/\varepsilon)&=H(x_0/\varepsilon)+
\frac{1}{\varepsilon}\int_0^t\int_0^{X^\varepsilon_v/\varepsilon}
\theta(s)dsdX^\varepsilon_v +\frac{1}{\varepsilon^{2(1-\kappa)}}
\int_0^t[\sigma^2(X^\varepsilon_s/\varepsilon)- \mathbf{a}]ds
\end{align*}
or, equivalently,
\begin{multline*}
\int_0^t[\sigma^2(X^\varepsilon_s/\varepsilon)-\mathbf{a}]ds =
\varepsilon^{2(1-\kappa)}\int_0^{X^\varepsilon_t/\varepsilon}\int_0^v
\theta(s)dsdv -
\varepsilon^{2(1-\kappa)}\int_0^{x_0/\varepsilon}\int_0^v
\theta(s)dsdv
-\\
\varepsilon^{1-2\kappa}\int_0^t\int_0^{X^\varepsilon_v/\varepsilon}
\theta(s) dsb(X^\varepsilon_v/\varepsilon) dv
-\varepsilon^{1-\kappa}\int_0^t\int_0^{X^\varepsilon_v/\varepsilon}
\theta(s)ds\sigma(X^\varepsilon_v/\varepsilon)dB_v.
\end{multline*}
Other steps of the proof repeat the previous ones and are omitted.
\qed

\section{MDP gap between oscillating and random environments}
\label{sec-6} The MDP regimes for oscillating and random
environments are proved for different ranges of $\kappa$:
$$
\kappa\in \Big(0,\frac{1}{2}\Big) \quad\text{and}\quad \kappa\in
\Big(0,\frac{1}{6}\Big)
$$
respectively. This fact is explained by faster homogenization effect
of the oscillating environment than the random one. This is clearly
seen from the  following proof sketch of the convergence in
\eqref{1f} for the oscillating environment.

Denote by $\theta(s)$ either $\frac{b(s)-\mathbf{b}}{\sigma^2(s)}$
or $1-\frac{\mathbf{a}}{\sigma^2(s)}$ and choose $\mathbf{b}$ and
$\mathbf{a}$ to satisfy the condition
\begin{equation}\label{zero}
\int_0^1\theta(s)ds=0,
\end{equation}
that is, $
\mathbf{b}=\int_0^1\frac{b(s)ds}{\sigma^2(s)}\big/\int_0^1\frac{ds}{\sigma^2(s)}
\quad\text{and}\quad
\mathbf{a}=1\big/\int_0^1\frac{1}{\sigma^2(s)}ds. $ Set $
H(x)=\int_0^x\int_0^v\theta(s)dsdv. $ Applying the It\^o formula to
$\varepsilon^{2(1-\kappa)}H(X^\varepsilon_t/\varepsilon)$, we find
that
\begin{align}\label{3.1a}
\left.
\begin{array}{lll}
\int_0^t[b(X^\varepsilon_s/\varepsilon)-\mathbf{b}]ds
\\
\int_0^t[\sigma^2(X^\varepsilon_s/\varepsilon)-\mathbf{a}]ds
\end{array}
\right\} &=\varepsilon^{2(1-\kappa)}
\int_0^{X^\varepsilon_t/\varepsilon}\int_0^v\theta(s)dsdv
\\
\label{3.2a} &\quad -\varepsilon^{2(1-\kappa)}
\int_0^{x_0/\varepsilon}\int_0^v\theta(s)dsdv
\\
\label{3.3a}&\quad
-\varepsilon^{1-2\kappa}\int_0^t\int_0^{X^\varepsilon_v/
\varepsilon}\theta(s)dsb(X^\varepsilon_v/\varepsilon)dv
\\
\nonumber &\quad
-\varepsilon^{1-\kappa}\int_0^t\int_0^{X^\varepsilon_v/
\varepsilon}\theta(s)ds\sigma (X^\varepsilon_v/\varepsilon)dB_v.
\end{align}
Since $\theta(s)$ is a periodic function, \eqref{zero} implies that
$\big|\int_0^t\theta(s)ds\big|$ is bounded uniformly in $t$. This is
the origin of strong homogenization. Namely, we have  the following
estimates for the terms in \eqref{3.1a} - \eqref{3.3a} (here $l$ is
a generic positive constant):
\begin{align*}
& \varepsilon^{2(1-\kappa)}
\Big|\int_0^{X^\varepsilon_t/\varepsilon}\int_0^v\theta(s) dsdv\Big|
\le l\varepsilon^{1-2\kappa}\sup_{t\le T}|X^\varepsilon_t|
\\
& \varepsilon^{2(1-\kappa)}
\Big|\int_0^{x_0/\varepsilon}\int_0^v\theta(s) dsdv\Big| \le
l\varepsilon^{1-2\kappa}|x_0|
\\
& \varepsilon^{1-2\kappa}\Big|\int_0^t\int_0^{X^\varepsilon_v/
\varepsilon}\theta(s)dsb(X^\varepsilon_v/\varepsilon)dv\Big|\le
l\varepsilon^{1-2\kappa}.
\end{align*}
The second and third upper bounds are deterministic and so the
corresponding terms are exponentially tight with the rate
$\varepsilon^{2\kappa}$ for any $\kappa>0$. Since $b$ and $\sigma$
are bounded, using \eqref{ccc}, we have
$$
\lim_{C\to\infty}\varlimsup_{\varepsilon\to
0}\varepsilon^{2\kappa}\log\mathsf{P}\big( \sup_{t\le
T}|X^\varepsilon_t|>C\big)=-\infty
$$
and, hence, the first term is exponentially tight for any $\kappa>0$
as well.

The restriction for $\kappa$ is imposed by the exponential
negligibility with the rate $\varepsilon^{2\kappa}$ of the
continuous martingale $
Z^\varepsilon_t:=\varepsilon^{1-\kappa}\sup_{t\le
T}\big|\int_0^t\int_0^{X^\varepsilon_v/
\varepsilon}\theta(s)ds\sigma
(X^\varepsilon_v/\varepsilon)dB_v\big|. $ Since
$$
\langle Z^\varepsilon\rangle_t=\varepsilon^{2(1-\kappa)}
\int_0^t\Big(\int_0^{X^\varepsilon_v/
\varepsilon}\theta(s)ds\Big)^2\sigma^2(X^\varepsilon_v/\varepsilon)dv
$$
and  $|\int_0^t\theta(s)ds|$ is uniformly bounded in $t$ (!),
$$
\langle Z^\varepsilon\rangle_T\le \varepsilon^{2(1-\kappa)}l.
$$
Therefore, by Lemma \ref{lem-2.1},
$$
\mathsf{P}\Big(\sup_{t\le T}\big|Z^\varepsilon_t\big|>\eta\Big)
=\mathsf{P}\Big(\sup_{t\le T}\big|Z^\varepsilon_t\big|>\eta, \langle
Z^\varepsilon \rangle_T \le \varepsilon^{2(1-\kappa)}l\Big) \le
2\exp\Big(-\frac{\eta^2} {2l\varepsilon^{2(1-\kappa)}}\Big),
$$
and so, $ \varepsilon^{2\kappa}\log \mathsf{P}\Big(\sup_{t\le
T}\big|Z^\varepsilon_t\big|> \eta\Big)\le \varepsilon^{2\kappa}\log
2 -\frac{\eta^2}
{2l\varepsilon^{2(1-2\kappa)}}\xrightarrow[\varepsilon\to
0]{}-\infty $ only if $\kappa<\frac{1}{2}$.

\appendix
\section{The LDP analysis}

The use of \eqref{1f} makes the proof of Theorem \ref{theo-000}
transparent. As was mentioned,  the implication
\begin{align}
\text{\eqref{1f}} \Rightarrow \text{the statement of Theorem
\ref{theo-000}} \label{A.11}
\end{align}
follows from Corollary 4.3.8  in \cite{puh2} (see also Corollary 6.7
in \cite{puh-ssr}) which are applicable not only to the setting
under consideration but also to various classes of semimartingales.
For reader's convenience, we show how \eqref{A.11} works in our
setting.

Let $X^\varepsilon=(X^\varepsilon_t)_{t\le T}$ be a continuous
semimartingale defined on a stochastic basis, with the general
conditions,
$(\varOmega,\mathbf{F},\mathcal{F}^\varepsilon=(\mathscr{F}^\varepsilon_t)_{t\le
T}, \mathsf{P})$:
\begin{equation*}
X^\varepsilon_t=x_0+\int_0^tb^\varepsilon_sds+\varepsilon^\kappa\int_0^t
\alpha^\varepsilon_sdB_s,
\end{equation*}
where the Brownian motion $B_t$  and the processes
$b^\varepsilon_t$, $\alpha^\varepsilon_t$ are
$\mathcal{F}^\varepsilon$-adapted (with
$\int_0^T|b^\varepsilon_t|dt<\infty$,
$\int_0^T(\alpha^\varepsilon_s)^2dt<\infty$), a.s.
$\varepsilon$ is a small positive parameter, $\kappa$ is a positive
number.

\begin{theorem}\label{theo-A.1'}
Assume $0<c_1\le|\alpha^\varepsilon_s|^2\le c_2$ and $|\beta^\varepsilon_s|\le c_3$ and there exist constants $\mathbf{b}$ and $\mathbf{a}>0$ such
that for any $\eta>0$
\begin{align}\label{gl1}
& \lim_{\varepsilon\to 0}\varepsilon^{2\kappa}\log\mathsf{P}
\Big(\sup_{t\le T}\Big|\int_0^T
[b^\varepsilon_s-\mathbf{b}]ds\Big|>\eta\Big)=-\infty
\\
& \lim_{\varepsilon\to 0}\varepsilon^{2\kappa}\log\mathsf{P}
\Big(\sup_{t\le T}\Big|\int_0^T
[(\alpha^\varepsilon_s)^2-\mathbf{a}]ds\Big|>\eta\Big)=-\infty.
\label{gl2}
\end{align}
Then, the family $ \{X^\varepsilon\}_{\varepsilon\to 0} $ obeys the
LDP with the rate $\varepsilon^{2\kappa}$ and the rate function
$$
J(u)=
  \begin{cases}
    \frac{1}{2\mathbf{a}}\int_0^T(\dot{u}_t-\mathbf{b})^2dt , & \begin{subarray}
    \phantom{}u_0=x_0\\du_t=\dot{u}_tdt,\\\int_0^T\dot{u}^2_tdt<\infty
    \end{subarray} \\
    \infty , & \text{otherwise}.
  \end{cases}
$$
\end{theorem}

The proof of this theorem uses a standard fact (see, e.g.
\cite{LipPuh}):
$$
\left.
  \begin{array}{ll}
  \text{Exponential Tightness} &
\\
\text{Local LDP} &
  \end{array}
\right\}\Rightarrow \ \text{LDP}.
$$

\medskip
\noindent For the proof of exponential tightness and  local LDP it
is convenient to use the stopping time
\begin{equation*}
\tau_{\varepsilon,\zeta}=\inf\Big\{t\le
T:\Big|\int_0^t[b^\varepsilon_s-\mathbf{b}]ds\Big|+
\Big|\int_0^t[(\alpha^\varepsilon_s)^2-\mathbf{a}]ds\Big|>\zeta\Big\}.
\end{equation*}
Notice also that \eqref{gl1} and \eqref{gl2} imply $
\lim_{\varepsilon\to
0}\varepsilon^{2\kappa}\log\mathsf{P}\big(\tau_{\varepsilon,\zeta}
<\infty\big)=-\infty $ and, therefore, for any measurable set
$\mathfrak{B}$,
\begin{equation}\label{bbb}
\varlimsup_{\varepsilon\to
0}\varepsilon^{2\kappa}\log\mathsf{P}\big(\mathfrak{B}\big)
\le\varlimsup_{\zeta\to 0}\varlimsup_{\varepsilon\to
0}\varepsilon^{2\kappa}\log\mathsf{P}\big(\mathfrak{B}
\cap\{\tau_{\varepsilon,\zeta}=\infty\}\big).
\end{equation}

\subsection{Exponential tightness}
\label{sec-A.1}

Following \cite{LipPuh}, we shall prove that
\begin{align}\label{expI}
&\lim_{C\to\infty}\varlimsup_{\varepsilon\to
0}\varepsilon^{2\kappa}\log\mathsf{P} \Big(\sup_{t\le
T}\big|X^\varepsilon_t\big|>C\Big)=-\infty
\\
& \lim_{\delta\to 0}\varlimsup_{\varepsilon\to
0}\varepsilon^{2\kappa}\log\sup_{\theta< T}\mathsf{P}\Big(
\sup_{0<t\le\delta}|X^\varepsilon_{\theta+t}-X^\varepsilon_\theta|>\eta\Big)=-\infty,
\forall \ \eta>0, \label{expII}
\end{align}
where $\theta$ is $\mathcal{F}^\varepsilon$-stopping time.

\subsubsection{\bf The proof of {\bf (\ref{expI})} }
\label{sec-A.1.1}

Set $\mathfrak{B}:=\big\{\sup_{t\le
T}\big|X^\varepsilon_t\big|>C\big\}$. Due to \eqref{bbb}, it
suffices to prove that $ \lim_{C\to\infty}\varlimsup_{\varepsilon\to
0}\varepsilon^{2\kappa}\log\mathsf{P} \big(\sup_{t\le
T}\big|X^\varepsilon_t\big|>C, \ \tau_{\varepsilon,\zeta}=\infty
\big)=-\infty. $ To this end, using the random variable $
|x_0|+|\mathbf{b}|T+\zeta+\sup_{t\le T}\big|\varepsilon^\kappa
\int_0^{t\wedge\tau_
{\varepsilon,\zeta}}\alpha^\varepsilon_sdB_s\big| $ as an upper
bound for $ \sup_{t\le T}|X^\varepsilon_t| $ on the set
$\{\tau_{\varepsilon,\zeta}=\infty\}$, the proof reduces to
$$
\lim_{C\to\infty}\varlimsup_{\varepsilon\to
0}\varepsilon^{2\kappa}\log\mathsf{P} \Big(\sup_{t\le
T}\big|\varepsilon^\kappa \int_0^{t\wedge\tau_
{\varepsilon,\zeta}}\alpha^\varepsilon_sdB_s\big|>C,\tau_{\varepsilon,\zeta}
=\infty\Big)=-\infty.
$$
Applying \eqref{ccc} to $ M_t=\varepsilon^\kappa
\int_0^{t\wedge\tau_ {\varepsilon,\zeta}}\alpha^\varepsilon_sdB_s $
and taking into account that
$$
\langle M\rangle_T=\varepsilon^{2\kappa} \int_0^{T\wedge\tau_
{\varepsilon,\zeta}}(\alpha^\varepsilon_s)^2ds\le
\varepsilon^{2\kappa}[\mathbf{a}T+\zeta],
$$
we find the following upper bound $ \mathsf{P}\Big(\sup_{t\le
T}\big|M_t\big|>C\Big) \le
2\exp\Big(-\frac{C^2}{2\varepsilon^{2\kappa}(\mathbf{a}T+\zeta)}\Big)
$
 which, in turn, provides \eqref{expI}.
\qed

\subsubsection{\bf The proof of {\bf (\ref{expII})}}
\label{sec-A.1.2}

Let $M_t$ be the same as above and let $
\mathcal{F}^{\varepsilon,\theta}=\big(\mathscr{F}
^\varepsilon_{\theta+t}\big)_{t\ge 0}. $ Denote
$N^\theta_t:=\big(M_{\theta+t}- M_\theta\big)$. The random process
$N^\theta_t$ is a martingale relative to $
\mathcal{F}^{\varepsilon,\theta} $ with $ \langle
N^\theta\rangle_t=\varepsilon^{2\kappa}\int_\theta^{\theta+t}
(\alpha^\varepsilon_s)^2ds. $ Denote $
L^\varepsilon_t:=\int_\theta^{\theta+t}b^\varepsilon_sds+N^\theta_t
$ and notice that \eqref{expII} is nothing but
$$
\lim_{\zeta\to 0}\varlimsup_{\delta\to 0}\varlimsup_{\varepsilon\to
0}
\varepsilon^{2\kappa}\log\sup_{\theta<T}\mathsf{P}\Big(\sup_{0<t\le\delta}
|L^\varepsilon_t|>\eta\Big)=-\infty,
$$
so that, in view of \eqref{bbb}, it suffices to prove that
$$
\lim_{\zeta\to 0}\varlimsup_{\delta\to 0}\varlimsup _{\varepsilon\to
0}
\varepsilon^{2\kappa}\log\sup_{\theta<T}\mathsf{P}\Big(\sup_{0<t\le\delta}
|L^\varepsilon_t|>\eta, \
\tau_{\varepsilon,\zeta}=\infty\Big)=-\infty.
$$
or, to verify the stronger condition
\begin{equation}\label{evv}
\lim_{\zeta\to 0}\varlimsup_{\delta\to 0}\varlimsup _{\varepsilon\to
0}
\varepsilon^{2\kappa}\log\sup_{\theta<T}\mathsf{P}\Big(\sup_{0<t\le\delta
\wedge\tau_{\varepsilon,\zeta}} |L^\varepsilon_t|>\eta\Big)=-\infty.
\end{equation}
Obviously, $ \sup_{0<t\le \delta}|L^\varepsilon_t|\le
\delta|\mathbf{b}|+\zeta+ \sup_{0<t\le
\delta}|N^\theta_{t\wedge\tau_{\varepsilon,\zeta}}|. $ For fixed
$\eta$ we choose sufficiently small $\delta$ and $\zeta$ such that
$\delta|\mathbf{b}|+\zeta<\eta$. Now, instead of \eqref{evv} we have
to prove
\begin{equation*}
\lim_{\zeta\to 0}\varlimsup_{\delta\to 0}\varlimsup _{\varepsilon\to
0}
\varepsilon^{2\kappa}\log\sup_{\theta<T}\mathsf{P}\Big(\sup_{0<t\le\delta
\wedge\tau_{\varepsilon,\zeta}}
|N^\theta_t|>\eta-\delta|\mathbf{b}|-\zeta\Big)=-\infty.
\end{equation*}
Further, due to $
 \langle N^\theta\rangle_{\delta\wedge\tau_{\varepsilon,\zeta}}
=\varepsilon^{2\kappa}\int_{\theta\wedge\tau_{\varepsilon,\zeta}}
^{(\theta+\delta\wedge\tau_{\varepsilon,\zeta})}
(\alpha^\varepsilon_s)^2ds\le\varepsilon^{2\kappa}[\delta\mathbf{a}+\zeta],
$ by applying \eqref{ccc} to $ N_{t\wedge\tau_
{\varepsilon,\zeta}}=\varepsilon^\kappa \int_0^{t\wedge\tau_
{\varepsilon,\zeta}}\alpha^\varepsilon_sdB_s $ we find the following
upper bound
$$
\mathsf{P}\Big(\sup_{t\le \delta\wedge\tau_ {\varepsilon,\zeta}}
\big|N^\theta_t\big|>\eta-\delta|\mathbf{b}|-\zeta\Big) \le
2\exp\Big(-\frac{(\eta-\delta|\mathbf{b}|-\zeta)^2}{2\varepsilon^{2\kappa}(\delta\mathbf{a}+\zeta)}\Big)
$$
which gives \eqref{expII}. \qed

\subsection{Local LDP}
\label{sec-A.2}

It is well known that for  exponentially tight family
$\{(X^\varepsilon_t)_{t\le T}\}_{\varepsilon\to 0}$ the rate
function coincides with the local rate function $J(u)$ determined by
the conditions: for any $u\in \mathbb{C}_{[0,T]}$,
\begin{align}\label{ub}
& \varlimsup_{\delta\to 0}\varlimsup_{\varepsilon\to
0}\varepsilon^{2\kappa}\log \mathsf{P}\Big(\sup_{t\le
T}|X^\varepsilon_t-u_t|\le\delta\Big)\le -J(u)
\\
& \varliminf_{\delta\to 0}\varliminf_{\varepsilon\to
0}\varepsilon^{2\kappa}\log \mathsf{P}\Big(\sup_{t\le
T}|X^\varepsilon_t-u_t|\le\delta\Big)\ge -J(u). \label{lb}
\end{align}

Since the range space of $J(u)$ is the interval $[0,\infty]$, we
compute separately $J(u)$ on two sets: $\mathfrak{U}_1=\{u\in
\mathbb{C}_{[0,T]}: J(u)<\infty\}$ and $\mathfrak{U}_2=\{u\in
\mathbb{C}_{[0,T]}: J(u)=\infty\}$.

\subsubsection{\bf Upper bound  under $\b{J(u)<\infty}$}\label{sec-A.2.2a}
\begin{lemma}\label{lem-A.1}
\eqref{ub} holds with
$J(u)=\frac{1}{2\mathbf{a}}\int_0^T\big(\dot{u}_t-\mathbf{b}
\big)^2dt$ for any
$$
u\in\Big\{u_0=x_0, \ du_t=\dot{u}_tdt, \
\int_0^T\dot{u}^2_tdt<\infty\Big\}=: \mathfrak{U}_1.
$$
\end{lemma}
\begin{proof}
Denote $\mathfrak{B}_\delta=\big\{\sup_{t\le
T}|X^\varepsilon_t-u_t|\le\delta\big\}$. In view of \eqref{bbb}, it
suffices to show that
\begin{align*}
\varlimsup_{\zeta\to 0}\varlimsup_{\delta\to
0}\varlimsup_{\varepsilon\to 0} \varepsilon^{2\kappa}
\log\mathsf{P}(\mathfrak{B}_\delta\cap\{\tau_{\varepsilon,\zeta}=\infty\})\le-
\frac{1}{2\mathbf{a}}\int_0^T(\dot{u}_t-\mathbf{b})^2dt.
\end{align*}
To this end, we introduce a martingale exponential
$$
\frak{z}_t=\exp\Big(\frac{1}{\varepsilon^{2\kappa}}\Big[\int_0^{T}
\lambda(s)[dX^\varepsilon_s-b^\varepsilon_sds]-\frac{1}{2}
\int_0^{T}\lambda^2(s)(a^\varepsilon_s)^2ds\Big]\Big), \ t\le T,
$$
where $\lambda(s)$ is a continuously differentiable function with
the derivative $\dot\lambda(s)$. Integrating by parts with the help
of It\^o's formula we find that $
\int_0^T\lambda(s)dX^\varepsilon_s=\lambda(T)X^\varepsilon_T-\lambda(0)x_0-
\int_0^TX^\varepsilon_s\dot{\lambda}(s)ds $ and rewrite
$\log\mathfrak{z}_t$ to the following form:
\begin{equation*}
\log\frak{z}_T
=\frac{1}{\varepsilon^{2\kappa}}\Big[\lambda(T)X^\varepsilon_T-
\lambda(0)x_0-\int_0^{T}X^\varepsilon_s\dot{\lambda}(s)ds
-\int_0^T\Big(\lambda(s)b^\varepsilon_s+\frac{\lambda^2(s)}{2}
(a^\varepsilon_s)^2\Big)ds\Big].
\end{equation*}
Write
\begin{align}\label{above}
& \lambda(T)X^\varepsilon_T-
\lambda(0)x_0-\int_0^{T}X^\varepsilon_s\dot{\lambda}(s)ds
-\int_0^T\Big(\lambda(s)b^\varepsilon_s+\frac{\lambda^2(s)}{2}
(a^\varepsilon_s)^2\Big)ds
\nonumber\\
&=\lambda(T)u_T- \lambda(0)u_0-\int_0^{T}u_s\dot{\lambda}(s)ds
-\int_0^T\Big(\lambda(s)\mathbf{b}+\frac{\lambda^2(s)}{2}
\mathbf{a}\Big)ds
\nonumber\\
&\quad
+\lambda(T)[X^\varepsilon_t-u_T]-\int_0^{T}[X^\varepsilon_s-u_s]\dot{\lambda}(s)ds
-\int_0^T\Big(\lambda(s)[b^\varepsilon_s-\mathbf{b}]+\frac{\lambda^2(s)}{2}
[(\alpha^\varepsilon_s)^2-\mathbf{a}]\Big)ds.
\end{align}
Integrating by parts we find that
\begin{align*}
&
\int_0^T\lambda(s)[b^\varepsilon_s-\mathbf{b}]ds=\lambda(T)\int_0^T[b^\varepsilon_s-
\mathbf{b}]ds-\int_0^T\dot{\lambda}(s)\int_0^s[b^\varepsilon_{s'}-\mathbf{b}]ds'ds
\\
& \int_0^T\frac{\lambda^2(s)}{2}[(a^\varepsilon_s)^2-\mathbf{a}]ds =
\frac{\lambda^2(T)}{2}\int_0^T[(a^\varepsilon_s)^2 -
\mathbf{a}]ds-\int_0^T\lambda(s)\dot{\lambda}(s)\int_0^s[(a^\varepsilon_{s'})^2
-\mathbf{a}]ds'ds,
\end{align*}
and transform \eqref{above} into
\begin{align*}
& \lambda(T)X^\varepsilon_T-
\lambda(0)x_0-\int_0^{T}X^\varepsilon_s\dot{\lambda}(s)ds
-\int_0^T\Big(\lambda(s)b^\varepsilon_s+\frac{\lambda^2(s)}{2}
(a^\varepsilon_s)^2\Big)ds
\\
&=\lambda(T)u_T- \lambda(0)u_0-\int_0^{T}u_s\dot{\lambda}(s)ds
-\int_0^T\Big(\lambda(s)\mathbf{b}+\frac{\lambda^2(s)}{2}
\mathbf{a}\Big)ds
\\
&\quad
+\lambda(T)[X^\varepsilon_t-u_T]-\int_0^{T}[X^\varepsilon_s-u_s]\dot{\lambda}(s)ds
\\
&\quad -\lambda(T)\int_0^T[b^\varepsilon_s-
\mathbf{b}]ds+\int_0^T\dot{\lambda}(s)\int_0^s[b^\varepsilon_{s'}-\mathbf{b}]ds'ds
\\
&\quad -\frac{\lambda^2(T)}{2}\int_0^T[(a^\varepsilon_s)^2 -
\mathbf{a}]ds+\int_0^T\lambda(s)\dot{\lambda}(s)\int_0^s[(a^\varepsilon_{s'})^2
-\mathbf{a}]ds'ds.
\end{align*}
The right hand side of this identity can be estimated from below on
the set $\mathfrak{B}_\delta\cap\{\tau_{\varepsilon,\zeta}=\infty\}$
by
\begin{align*}
&\lambda(T)u_T- \lambda(0)u_0-\int_0^{T}u_s\dot{\lambda}(s)ds
-\int_0^T\Big(\lambda(s)\mathbf{b}+\frac{\lambda^2(s)}{2}
\mathbf{a}\Big)ds
\\
&\quad -\delta\Big(|\lambda(T)|+\int_0^{T}|\dot{\lambda}(s)|ds\Big)
\\
&\quad -\zeta\Big(|\lambda(T)|+\int_0^T|\dot{\lambda}(s)|ds
+\frac{\lambda^2(T)}{2}+\int_0^T|\lambda(s)\dot{\lambda}(s)|ds\Big).
\end{align*}
With $ l_1=|\lambda(T)|+\int_0^{T}|\dot{\lambda}(s)|ds+
|\lambda(T)|+\int_0^T|\dot{\lambda}(s)|ds
+\frac{\lambda^2(T)}{2}+\int_0^T|\lambda(s)\dot{\lambda}(s)|ds $
and the identity $ \lambda(T)u_T-
\lambda(0)u_0-\int_0^{T}u_s\dot{\lambda}(s)ds=\int_0^T\lambda(s)\dot{u}_tdt
$ we find the following lower bound for $\mathfrak{z}_T$:
\begin{equation}\label{z+z}
\mathfrak{z}_*=\exp\Big(\frac{1}{\varepsilon^{2\kappa}}\Big[\int_0^T\Big[\lambda(s)(\dot{u}_s-\mathbf{b})-
\frac{\lambda^2(s)}{2}\Big]ds-l_1(\delta+\zeta)\Big]\Big).
\end{equation}
The martingale
exponential $\mathfrak{z}_t$ is a positive local martingale and a
supermartingale too. Hence, $\mathsf{E}\mathfrak{z}_T\le 1$. This
bound implies $ \mathsf{E}I_{\{\mathfrak{B}_\delta
\cap(\tau_{\varepsilon,\zeta}=\infty)\}}\mathfrak{z}_T\le 1 $ and,
in turn,
$$
\mathsf{E}I_{\{\mathfrak{B}_\delta
\cap(\tau_{\varepsilon,\zeta}=\infty)\}}\mathfrak{z}_*\le 1.
$$
Jointly with \eqref{z+z} the latter implies
\begin{equation}\label{nad}
\varlimsup_{\zeta\to 0}\varlimsup_{\delta\to
0}\varlimsup_{\varepsilon\to 0} \varepsilon^{2\kappa}
\log\mathsf{P}(\mathfrak{B}_\delta
\cap\{\tau_{\varepsilon,\zeta}=\infty\})\le-
\int_0^{T}\Big\{\lambda(s)[\dot{u}_s
-\mathbf{b}]-\frac{\lambda^2(s)}{2} \mathbf{a}\Big\}ds.
\end{equation}
Recall that \eqref{nad} is valid provided that $\lambda(s)$ is
a continuously differentiable function.
Assume that  $\dot{u}_s$ is also  continuously differentiable.
 Then taking
$ \lambda(s)\equiv\frac{\dot{u}_s-\mathbf{b}}{\mathbf{a}} $ we
obtain the desired result. In the general case, $\dot{u}_t$ is  only
a density of $u_t$ relative to $dt$, so that, $\lambda(t)$ as chosen
above may not be continuously differentiable. In this case we use
the identity
$$
-\int_0^T\Big\{\lambda_m(s)[\dot{u}_s
-\mathbf{b}]-\frac{\lambda^2_m(s)}{2}
\mathbf{a}\Big\}ds=-\int_0^T\frac{(\dot{u}_s-\mathbf{b})^2}{2\mathbf{a}}ds+
\int_0^T\frac{\mathbf{a}}{2}
\Big(\lambda_m(s)-\frac{\dot{u}_s-\mathbf{b}}{\mathbf{a}}\Big)^2ds,
$$
where $\lambda_m(s)$ is a sequence of continuously differentiable
functions such that
$$
\lim_{m\to\infty}\int_0^T\frac{\mathbf{a}}{2}
\Big(\lambda_m(s)-\frac{\dot{u}_s-\mathbf{b}}{\mathbf{a}}\Big)^2ds=0.
$$
\end{proof}

\subsubsection{\bf Upper bound  under $\b{J(u)=\infty}$}

Since $ \mathfrak{U}_2=\mathbb{C}_{[0,T]}\setminus\mathfrak{U}_1, $
it suffices to verify the upper bound in \eqref{ub} under the
following conditions:
\begin{enumerate}
\item[{\bf (c.1)}] $u_0=x_0, \ du_t\ll dt, \ \int_0^T\dot{u}^2dt=\infty$
\item[{\bf (c.2)}] $u_0=x_0, \ du_t\not\ll dt$
\item[{\bf (c.3)}] $u_0\ne x_0$
\end{enumerate}
\begin{lemma}\label{lem-A.2}
For any of {\bf (c.1)}, {\bf (c.2)} and {\bf (c.3)}, the upper bound
in \eqref{ub} holds with $J(u)=\infty$.
\end{lemma}

\begin{proof}

\mbox{}

{\bf (c.1)} By \eqref{nad} and \eqref{bbb}, for any continuously
differentiable function $\lambda(s)$, we have
\begin{equation}\label{nad1}
\varlimsup_{\delta\to 0}\varlimsup_{\varepsilon\to 0}
\varepsilon^{2\kappa}
\log\mathsf{P}\big(\mathfrak{B}_\delta\big)\le-
\int_0^{T}\Big\{\lambda(s)[\dot{u}_s
-\mathbf{b}]-\frac{\lambda^2(s)}{2} \mathbf{a}\Big\}ds.
\end{equation}

Let us take $\lambda_n(s)=\frac{\dot{u}_s-\mathbf{b}}{\mathbf{a}}
I_{\{|\dot{u}_s|\le n\}}$ and choose a sequence of continuously
differentiable functions $\lambda_{m,n}(s)$ such that $
\lim_{m\to\infty}\int_0^T[\lambda_n(s)-\lambda_{m,n}(s)]^2=0. $ For
$\lambda_{m,n}(s)$ \eqref{nad1} implies:
\begin{equation}\label{nad2}
\varlimsup_{\delta\to 0}\varlimsup_{\varepsilon\to 0}
\varepsilon^{2\kappa}
\log\mathsf{P}\big(\mathfrak{B}_\delta\big)\le-
\int_0^{T}\Big\{\lambda_{m,n}(s)[\dot{u}_s
-\mathbf{b}]-\frac{\lambda^2_{m,n}(s)}{2} \mathbf{a}\Big\}ds.
\end{equation}
The right hand side of \eqref{nad2} converges to $ -
\int_0^{T}\big\{\lambda_{n}(s)[\dot{u}_s
-\mathbf{b}]-\frac{\lambda^2_{n}(s)}{2} \mathbf{a}\big\}ds $ with
$m\to\infty$. Noticing that $ -
\int_0^{T}\big\{\lambda_{n}(s)[\dot{u}_s
-\mathbf{b}]-\frac{\lambda^2_{n}(s)}{2} \mathbf{a}\big\}ds= -
\frac{1}{2\mathbf{a}}\int_0^{T}\big(\dot{u}_s-\mathbf{b}\big)^2I_{\{|\dot{u}_s|\le
n\}} ds, $ we find that for any $n\ge 1$,
$$
\varlimsup_{\delta\to 0}\varlimsup_{\varepsilon\to 0}
\varepsilon^{2\kappa}
\log\mathsf{P}\big(\mathfrak{B}_\delta\big)\le-
\frac{1}{2\mathbf{a}}\int_0^{T}\big(\dot{u}_s-\mathbf{b}\big)^2I_{\{|\dot{u}_s|\le
n\}} ds\xrightarrow[n\to\infty]{}-\infty.
$$

\medskip
{\bf (c.2)} We show that $du_t\not\ll dt$ enables us to choose
a sequence $[s_i,t_i)$, $i=1,\ldots n$ of nonoverlapping intervals on $[0,T]$
such that
\begin{equation}\label{ee}
\lim_{n\to\infty}
\sum_{i=1}^n\frac{|u_{t^n_i}-u_{s^n_i}|^2}{t^n_i-s^n_{i}}=\infty.
\end{equation}
By
the Cauchy-Schwarz inequality
\begin{equation}\label{CS}
\Big(\sum\limits_{i=1}^n|u_{t^n_i}-u_{s^n_{i}}|\Big)^2\le
\sum_{i=1}^n[t^n_i-s^n_i]
\sum_{i=1}^n\frac{|u_{t^n_i}-u_{s^n_{i}}|^2}{t^n_i-s^n_{i}}.
\end{equation}
Further, for any small positive number $\gamma$ one can choose
intervals $[s_i,t_i)$, $i=1,\ldots n$ such that $
\sum_{i=1}^n[t^n_i-s^n_i]\le \gamma $ for any $n\ge n_\gamma$, where
$n_\gamma$ is some number depending on $\gamma$ and, at the same
time, $ \sum_{i=1}^n|u_{t^n_i}-u_{s^n_i}|\ge D>0. $ Hence,
\eqref{CS} implies $ \frac{D^2}{\gamma}\le
\sum_{i=1}^n\frac{|u_{t^n_i}-u_{s^n_{i}}|^2}{t^n_i-s^n_{i}} $
and \eqref{ee} holds in view of $D$ remains strictly positive with $\gamma\to 0$.

We prove that
\begin{equation*}
\varlimsup_{\zeta\to 0}\varlimsup_{\delta\to
0}\varlimsup_{\varepsilon\to 0}\log
\mathsf{P}\big(\mathfrak{B}_\delta\cap\{\tau_{\varepsilon,\zeta}
=\infty\}\big) \le
-\frac{1}{2\mathbf{a}}\Big\{\mathbf{b}^2T-2\mathbf{b}(u_T-u_0)+
 \sum_{i=1}^n\frac{[u_{t^n_i}-u_{s^n_{i}}]^2}{t^n_i-s^n_{i}}\Big\}
\end{equation*}
and, then, apply \eqref{ee}.

With $ \lambda(t)=\sum_i\lambda_iI_{\{s^n_i\le t< t^n_{i}\}}, $
we introduce a martingale exponential
$$
\frak{z}_t=\exp\Big(\frac{1}{\varepsilon^{2\kappa}}\Big[\int_0^{T}
\lambda(s)[dX^\varepsilon_s-b^\varepsilon_sds]-\frac{1}{2}
\int_0^{T}\lambda^2(s)(a^\varepsilon_s)^2ds\Big]\Big)
$$
and estimate  it from below on the set
$\mathfrak{B}_\delta\cap\{\tau_{\varepsilon, \zeta}=\infty\}$. Write
\begin{equation*}
\begin{aligned}
\log\frak{z}_T&
=\frac{1}{\varepsilon^{2\kappa}}\sum_{i=1}^n\Big[\lambda_i[X^\varepsilon_{t^n_i}-
X^\varepsilon_{s^n_i}]-\int_{s^n_i}^{t^n_i}\Big\{\lambda_ib^\varepsilon_s+
\frac{1}{2}\lambda^2_i(a^\varepsilon_s)^2\Big\}ds\Big]
\\
&\ge
\frac{1}{\varepsilon^{2\kappa}}\sum_{i=1}^n\Big\{\lambda_i[u_{t^n_i}-u_{s^n_i}]
-\Big[\lambda_i\mathbf{b}+\frac{\lambda^2_i}{2}\mathbf{a}\Big](t^n_i-s^n_i)\Big\}
\\
& \quad
-\frac{1}{\varepsilon^{2\kappa}}\sum_{i=1}^n\Big(2|\lambda_i|\delta+\zeta\Big[|\lambda_i|
+\frac{\lambda^2_i}{2}\Big]
\Big)(t^n_i-s^n_i):=\log\mathfrak{z}_*.
\end{aligned}
\end{equation*}
Since $ \mathsf{E}\mathfrak{z}_T\le 1,
$ also $ \mathsf{E}I_{\{\mathfrak{B}_\delta\cap\{\tau_{\varepsilon,
\zeta}=\infty\}\}}\mathfrak{z}_T\le 1 $ and, in turn, $
\mathsf{E}I_{\{\mathfrak{B}_\delta\cap\{\tau_{\varepsilon,
\zeta}=\infty\}\}}\mathfrak{z}_*\le 1. $ The latter implies
\begin{multline*}
\varlimsup_{\zeta\to 0}\varlimsup_{\delta\to
0}\varlimsup_{\varepsilon\to 0}
\varepsilon^{2\kappa}\log\mathsf{P}(\mathfrak{B}_\delta\cap\{\tau_{\varepsilon,\zeta=\infty}\})
\\
\le- \sum_{i=1}^n\Big[\lambda_i([u_{t^n_i}-
u_{s^n_i}]-\mathbf{b}[t^n_i-s^n_i])-
\frac{\lambda^2_i}{2}\mathbf{a}[t^n_i-s^n_i]\Big]
\end{multline*}
and it is left to choose $
 \lambda_i=\frac{[u_{t^n_i}-
u_{s^n_i}]-\mathbf{b}[t^n_i-s^n_i]}{\mathbf{a}[t^n_i-s^n_i]}
$ and apply \eqref{ee}.

\medskip
{\bf (c.3)} is obvious.
\end{proof}

\subsubsection{\bf Lower bound}
\label{sec-A.2.4}

Since the upper bound equals $-\infty$ for any $u\in\mathfrak{U}_2$,
the lower bound has to be checked for $u\in\mathfrak{U}_1$ only.
\begin{lemma}
\eqref{lb} holds with
$J(u)=\frac{1}{2\mathbf{a}}\int_0^T\big(\dot{u}_t-\mathbf{b}
\big)^2dt$ for any $ u\in\mathfrak{U}_1. $
\end{lemma}
\begin{proof}
It suffices to prove
\begin{equation}\label{LB}
\varliminf_{\delta\to 0}\varliminf_{\varepsilon\to 0}
\varepsilon^{2\kappa}\log
\mathsf{P}\big(\mathfrak{B}_\delta\big)\ge- \frac{1}{2\mathbf{a}}
\int_0^T[\dot{u}_t-\mathbf{b}]^2dt
\end{equation}
only for $u$ with continuous second derivative $\ddot{u}_t$.
 Indeed,
as $\int_0^T\dot{u}^2_tdt<\infty$, there exists a sequence $
\big(u_m(t)\big)_{m\ge 1} $ of twice continuously differentiable
functions such that their derivatives $\big(\dot{u}_m(t)\big)_{m\ge
1}$ approximate $\dot{u}_t$: $
\lim_{m\to\infty}\int_0^T\big(\dot{u}_t-\dot{u}_m(t)\big)^2dt=0. $
The latter also provides
\begin{equation}\label{nepovtor}
 \sup_{t\le T}|u_t-u_m(t)|\le
\big(T\int_0^T[\dot{u}_s-\dot{u}_m(s)]^2ds
\big)^{1/2}\xrightarrow[m\to\infty]{}0.
\end{equation}
Set $ \mathfrak{B}_{\delta,m}=\big\{\sup_{t\le
T}|X^\varepsilon_t-u_m(t)|\le \delta\big\} $ and
suppose we know that
\begin{equation*}
\varliminf_{\delta\to 0}\varliminf_{\varepsilon\to 0}
\varepsilon^{2\kappa}\log
\mathsf{P}\big(\mathfrak{B}_{\delta,m}\big)\ge-
\frac{1}{2\mathbf{a}} \int_0^T[\dot{u}_m(t)-\mathbf{b}]^2dt.
\end{equation*}
For $\gamma>\delta>0$,
 write
\begin{align*}
\mathfrak{B}_{\delta,m} & \subseteq \Big(\mathfrak{B}_{\delta,m}\cap
\Big\{\sup_{t\le T}\big|u_t-u_m(t)\big|\le\gamma\Big\}\Big) \bigcup
\Big\{\sup_{t\le T}\big|u_t-u_m(t)\big|>\gamma\Big\}
\\
&\subseteq \Big\{\sup_{t\le T}\Big|X^\varepsilon_t-u_t\Big|\le
2\gamma\Big\} \bigcup \Big\{\sup_{t\le
T}\big|u_t-u_m(t)\big|>\gamma\Big\}.
\end{align*}
For sufficiently large $m$, the
set $ \big\{\sup_{t\le T}|u_t-u_m(t)|>\gamma\big\} $ is empty.

Hence, $ \mathsf{P}\big(\sup_{t\le
T}\big|X^\varepsilon_t-u_t\big|\le 2\gamma\big) \ge
\mathsf{P}\big(\mathfrak{B}_{\delta,m}\big) $ and, therefore,
$$
\varliminf_{\varepsilon\to 0}\varepsilon^{2\kappa}\log
\mathsf{P}\Big(\sup_{t\le T}\Big|X^\varepsilon_t-u_t\Big|\le 2\gamma
\Big)\ge- \frac{1}{2\mathbf{a}} \int_0^T[\dot{u}_t-\mathbf{b}]^2dt.
$$
and it is left to pass to the limit with $\gamma\to 0$.

\medskip
The second helpful fact is that \eqref{LB} with $u$, having
continuously differentiable $\ddot{u}$,  follows from
\begin{equation}\label{LBz}
\varliminf_{\zeta\to 0}\varliminf_{\delta\to \infty}
\varliminf_{\varepsilon\to 0} \varepsilon^{2\kappa}\log
\mathsf{P}\Big(\sup_{t\le
T\wedge\tau_{\varepsilon,\zeta}}|X^\varepsilon_t- u_t|\le
\delta\Big)\ge- \frac{1}{2\mathbf{a}}
\int_0^T[\dot{u}_t-\mathbf{b}]^2dt.
\end{equation}
Write
\begin{align*}
\Big\{\sup_{t\le T\wedge\tau_{\varepsilon,\zeta}}|X^\varepsilon_t-
u_t|\le \delta\Big\}&=\Big\{\sup_{t\le T}|X^\varepsilon_t- u_t|\le
\delta,\tau_{\varepsilon,\zeta}=\infty\big\}
\\
&\quad \bigcup \Big\{\sup_{t\le
T\wedge\tau_{\varepsilon,\zeta}}|X^\varepsilon_t- u_t|\le
\delta,\tau_{\varepsilon,\zeta}<\infty\Big\}
\\
& \subseteq \Big\{\sup_{t\le T}|X^\varepsilon_t- u_t|\le
\delta\Big\} \bigcup \big\{\tau_{\varepsilon,\zeta}<\infty\big\}.
\end{align*}
Hence, we obtain the following inequality
\begin{align*}
\mathsf{P}\Big(\sup_{t\le
T\wedge\tau_{\varepsilon,\zeta}}|X^\varepsilon_t- u_t|\le
\delta\Big) \le 2\Big[\mathsf{P}\Big(\sup_{t\le T}|X^\varepsilon_t-
u_t|\le
\delta\Big)\bigvee\mathsf{P}\big(\tau_{\varepsilon,\zeta}<\infty\big)\Big]
\end{align*}
and in view of  \eqref{LBz} and
$
\varlimsup_{\zeta\to\infty}\varlimsup_{\varepsilon\to
0}\varepsilon^{2\kappa} \log
\mathsf{P}\big(\tau_{\varepsilon,\zeta}<\infty\big)=-\infty$,
we obtain the required \eqref{LB}.

\medskip
Henceforth, we focus on the proof of \eqref{LBz}, where  $u_t$ is
assumed to be twice continuously differentiable function.

Let us introduce a martingale exponential $
\mathfrak{z}_t=e^{M_t-\frac{1}{2}\langle M\rangle_t} $ with
$$
M_t=\frac{1}{\varepsilon^{\kappa}}\int_0^{t\wedge\tau_{\varepsilon,\zeta}}
\frac{\dot{u}_t-\beta^\varepsilon_s}{\alpha^\varepsilon_s}dB_s
\quad \text{and} \quad \langle M\rangle_t=
\frac{1}{\varepsilon^{2\kappa}}
\int_0^{t\wedge\tau_{\varepsilon,\zeta}}
\frac{(\dot{u}_t-\beta^\varepsilon_s)^2}{(\alpha^\varepsilon_t)^2}dt.
$$
Since $\dot{u}_t$ is bounded, $\langle M\rangle_T$ is bounded too. Hence, we have
$ \mathsf{E}\mathfrak{z}_T=1. $ Set $
d\bar{\mathsf{P}}=\mathfrak{z}_Td\mathsf{P}. $ Owing to
$\mathfrak{z}_T>0$, $\mathsf{P}$- a.s., we have $
\bar{\mathsf{P}}\sim\mathsf{P} $ with $
d\mathsf{P}=\mathfrak{z}^{-1}_Td\bar{\mathsf{P}}. $ Write
\begin{equation}\label{A.191}
\mathsf{P}\Big(\sup_{t\le
T\wedge\tau_{\varepsilon,\zeta}}|X^\varepsilon_t- u_t|\le
\delta\Big)= \int\limits_{\big\{\sup_{t\le
T\wedge\tau_{\varepsilon,\zeta}}|X^\varepsilon_t -u_t|\le
\delta\big\}} \mathfrak{z}^{-1}_Td\bar{\mathsf{P}}.
\end{equation}
First we show  that
\begin{equation}\label{fst}
\varliminf_{\zeta\to 0}\varliminf_{\varepsilon\to 0}
\bar{\mathsf{P}}\Big(\sup_{t\le
T\wedge\tau_{\varepsilon,\zeta}}|X^\varepsilon_t -u_t|\le
\delta\Big)=1.
\end{equation}
By the Girsanov theorem (see, e.g., Theorem 2, \S 5, Ch. 4 in
\cite{LSMar}), the random process
\begin{equation}\label{A.201}
\bar{B}_t=B_t-\frac{1}{\varepsilon^\kappa}\int_0^{t\wedge\tau_{\varepsilon,\zeta}}
\frac{\dot{u}_s-\beta^\varepsilon_s}{\alpha^\varepsilon_s}ds
\end{equation}
is a Brownian motion under $\bar{\mathsf{P}}$. Both process
$X^\varepsilon_t$ and $\mathfrak{z}_t$ are semimartingales under
$\bar{\mathsf{P}}$. In particular,
\begin{align*}
 X^\varepsilon_{t\wedge\tau_{\varepsilon,\zeta}}=x_0+
\int_0^{t\wedge\tau_{\varepsilon,\zeta}}\dot{u}_sds+\varepsilon^\kappa
\int_0^{t\wedge\tau_{\varepsilon,\zeta}}\alpha^\varepsilon_sd\bar{B}_s
=u_{t\wedge\tau_{\varepsilon,\zeta}}+\varepsilon^\kappa
\int_0^{t\wedge\tau_{\varepsilon,\zeta}}\alpha^\varepsilon_sd\bar{B}_s.
\end{align*}
Consequently,
$$
\sup_{t\le
T\wedge\tau_{\varepsilon,\zeta}}\big|X^\varepsilon_t-u_t\big|=
\sup_{t\le T}\Big|\varepsilon^\kappa
\int_0^{t\wedge\tau_{\varepsilon,\zeta}}\alpha^\varepsilon_sd\bar{B}_s\Big|.
$$
Using the Doob inequality, we find that
\begin{align}
& \bar{\mathsf{P}}\Big(\sup_{t\le
T\wedge\tau_{\varepsilon,\zeta}}|X^\varepsilon_t -u_t|\le\delta\Big)=
\bar{\mathsf{P}}\Big(
\sup_{t\le
T}\Big|\varepsilon^\kappa
\int_0^{t\wedge\tau_{\varepsilon,\zeta}}\alpha^\varepsilon_sd\bar{B}_s\Big|\le\delta
\Big)
\nonumber\\
&=1-\bar{\mathsf{P}}\Big(\sup_{t\le T}\Big|\varepsilon^\kappa
\int_0^{t\wedge\tau_{\varepsilon,\zeta}}\alpha^\varepsilon_sd\bar{B}_s\Big|
>\delta\Big)
\nonumber\\
&\ge 1-\frac{4\varepsilon^{2\kappa}}{\delta^2}
\bar{\mathsf{E}}\int_0^{T\wedge\tau_{\varepsilon,\zeta}}(\alpha^\varepsilon_s)^2ds
\ge 1-\frac{4\varepsilon^{2\kappa}}{\delta^2}
[\mathbf{a}T+\zeta]\xrightarrow[\varepsilon\to 0]{}1. \label{A.223}
\end{align}

\medskip
We return to the proof of \eqref{LBz}.
Due to \eqref{A.201}, we
have
\begin{align*}
\mathfrak{z}_T=\exp\Big(\frac{1}{\varepsilon^{\kappa}}
\int_0^{T\wedge\tau_{\varepsilon,
\zeta}}\frac{\dot{u}_t-\beta^\varepsilon_s}{\alpha^\varepsilon_s}d\bar{B}_s
+\frac{1}{2\varepsilon^{2\kappa}}\int_0^{T\wedge\tau_{\varepsilon,
\zeta}}\frac{[\dot{u}_t-\beta^\varepsilon_s]^2}{(\alpha^\varepsilon_s)^2}ds\Big)
\end{align*}
and, setting
 $
\mathfrak{B}_{\delta,\varepsilon}:=\{\sup_{t\le
T\wedge\tau_{\varepsilon,\zeta}}|X^\varepsilon_t- u_t|\le \delta\},
$
transform \eqref{A.191} to
\begin{gather*}
 \mathsf{P}\big(\mathfrak{B}_{\delta,\varepsilon}\big)
= \int_{\mathfrak{B}_{\delta,\varepsilon}}
\exp\Big(-\frac{1}{\varepsilon^{2\kappa}}\int_0^{T\wedge\tau_{\varepsilon,
\zeta}}\frac{\dot{u}_t-\beta^\varepsilon_s}{\alpha^\varepsilon_s}d\bar{B}_s
-\frac{1}{2\varepsilon^{2\kappa}}\int_0^{T\wedge\tau_{\varepsilon,
\zeta}}\frac{[\dot{u}_s-\beta^\varepsilon_s]^2}{(\alpha^\varepsilon_s)^2}ds\Big)
d\bar{\mathsf{P}}.
\end{gather*}
Further, it is convenient to apply a few obvious relations.
We choose $l$ as un upper bound for
$$
\frac{(\dot{u}_s-\mathbf{b})^2}
{\mathbf{a}(\alpha^\varepsilon_s)^2}
+\frac{2|\dot{u}_s-\mathbf{b}|+
|\mathbf{b}-\beta^\varepsilon_s|}
{\mathbf{a}}\Big[1+\frac{|\mathbf{a}-(\alpha^\varepsilon_s)^2|}{(\alpha^\varepsilon_s)^2}
\Big].
$$
Write
\begin{align*}
\int_0^{T\wedge\tau_{\varepsilon,\zeta}}\frac{[\dot{u}_s-\beta^\varepsilon_s]^2}
{(\alpha^\varepsilon_s)^2}ds
&=
\int_0^{T\wedge\tau_{\varepsilon,\zeta}}
\frac{[(\dot{u}_s-\mathbf{b})+(\mathbf{b}-\beta^\varepsilon_s)]^2}
{\mathbf{a}}\Big[1+\frac{\mathbf{a}-(\alpha^\varepsilon_s)^2}{(\alpha^\varepsilon_s)^2}
\Big]ds
\\
&\le \int_0^{T\wedge\tau_{\varepsilon,\zeta}}
\frac{(\dot{u}_s-\mathbf{b)^2}}{\mathbf{a}}ds
+\int_0^{T\wedge\tau_{\varepsilon,\zeta}}
\frac{(\dot{u}_s-\mathbf{b})^2}
{\mathbf{a}}\frac{|\mathbf{a}-(\alpha^\varepsilon_s)^2|}{(\alpha^\varepsilon_s)^2}
ds
\\
&\quad
+\int_0^{T\wedge\tau_{\varepsilon,\zeta}}
\frac{2|\dot{u}_s-\mathbf{b}||\mathbf{b}-\beta^\varepsilon_s|+
(\mathbf{b}-\beta^\varepsilon_s)^2}
{\mathbf{a}}\Big[1+\frac{|\mathbf{a}-(\alpha^\varepsilon_s)^2|}{(\alpha^\varepsilon_s)^2}
\Big]ds
\\
&\le \int_0^{T\wedge\tau_{\varepsilon,\zeta}}
\frac{(\dot{u}_s-\mathbf{b)^2}}{\mathbf{a}}ds
+l\int_0^{T\wedge\tau_{\varepsilon,\zeta}}\Big[|\mathbf{a}-(\alpha^\varepsilon_s)^2|+
|\mathbf{b}-\beta^\varepsilon_s|\Big]ds
\\
&\le \int_0^{T\wedge\tau_{\varepsilon,\zeta}}
\frac{(\dot{u}_s-\mathbf{b)^2}}{\mathbf{a}}ds
+\zeta l.
\end{align*}
The latter implies
\begin{align*}
&
\int_{\mathfrak{B}_{\delta,\varepsilon}}
\exp\Big(-\frac{1}{\varepsilon^{\kappa}}\int_0^{T\wedge\tau_{\varepsilon,
\zeta}}\frac{\dot{u}_s-\beta^\varepsilon_s}{\alpha^\varepsilon_s}d\bar{B}_s
-\frac{1}{2\varepsilon^{2\kappa}}\int_0^{T\wedge\tau_{\varepsilon,
\zeta}}\frac{[\dot{u}_s-\beta^\varepsilon_s]^2}{(\alpha^\varepsilon_s)^2}ds\Big)
d\bar{\mathsf{P}}
\\
&\ge
\exp\Big(-\frac{1}{2\varepsilon^{2\kappa}}\int_0^T
\frac{(\dot{u}_s-\mathbf{b})^2}{\mathbf{a}}ds-\frac{\zeta l}{2\varepsilon^{2\kappa}}\Big)
\int_{\mathfrak{B}_{\delta,\zeta}}1\wedge
\exp\Big(-\frac{1}{\varepsilon^\kappa}\int_0^{T\wedge\tau_{\varepsilon,
\zeta}}\frac{\dot{u}_t-\beta^\varepsilon_s}{(\alpha^\varepsilon_s)^2}d\bar{B}\Big)
d\bar{\mathsf{P}}
\end{align*}
(``$1\wedge\exp(\cdots)$'' is introduced in order to have bounded integrand in the last
integral above).

Thus,
 we obtain the following lower bound:
\begin{multline*}
 \mathsf{P}^{\varepsilon^{2\kappa}}\Big(\mathfrak{B}_{\delta,\varepsilon}\Big)
\ge\exp\Big(-\frac{1}{2\mathbf{a}}\int_0^T
(\dot{u}_s-\mathbf{b})^2dt-\frac{\zeta l}{2}\Big)
\\
\times\Bigg( \int_{\mathfrak{B}_{\delta,\varepsilon}}1\wedge
\exp\Big(-\frac{1}{\varepsilon^\kappa}\int_0^{T\wedge\tau_{\varepsilon,
\zeta}}\frac{\dot{u}_s-\beta^\varepsilon_s}{(\alpha^\varepsilon_s)^2}d\bar{B}\Big)
d\bar{\mathsf{P}}\Bigg)^{\varepsilon^{2\kappa}}.
\end{multline*}
Further, by the H\"older inequality (here $\varepsilon<1$)
\begin{align*}
&\int\limits_{\mathfrak{B}_{\delta,\varepsilon}}1\wedge
\exp\Big(-\varepsilon^\kappa\int_0^{T\wedge\tau_{\varepsilon,
\zeta}}\frac{\dot{u}_s-\beta^\varepsilon_s}{(\alpha^\varepsilon_s)^2}d\bar{B}_s
\Big)d\bar{\mathsf{P}}
\\
&\le \Bigg(\int\limits_{\mathfrak{B}_{\delta,\varepsilon}}
1\wedge\exp\Big(-\frac{1}{\varepsilon^{\kappa}}\int_0^{T\wedge\tau_{\varepsilon,
\zeta}}\frac{\dot{u}_s-\beta^\varepsilon_s}{(\alpha^\varepsilon_s)^2}d\bar{B}_s
\Big)d\bar{\mathsf{P}}\Bigg)^{\varepsilon^{2\kappa}}.
\end{align*}

Thus, we obtain
\begin{align*}
& \varliminf_{\varepsilon\to 0}\varepsilon^{2\kappa}\log\mathsf{P}
\Big(\sup_{t\le T\wedge\tau_{\varepsilon,\zeta}}|X^\varepsilon_t-
u_t|\le \delta\Big)\ge-\frac{1}{2\mathbf{a}}\int_0^T
(\dot{u}_s-\mathbf{b})^2ds
\\
&\quad -\frac{\zeta l}{2}+\varliminf_{\varepsilon\to
0}\varepsilon^{2\kappa}\log \int_{\mathfrak{B}_{\delta,\varepsilon}}
1\wedge\exp\Big(-\varepsilon^\kappa\int_0^{T\wedge\tau_{\varepsilon,
\zeta}}\frac{\dot{u}_s-\beta^\varepsilon_s}{(\alpha^\varepsilon_s)^2}d\bar{B}_s
\Big)d\bar{\mathsf{P}}
\\
&\ge-\frac{1}{2\mathbf{a}}\int_0^T
(\dot{u}_s-\mathbf{b})^2ds-\frac{\zeta l}{2}+\log
\varliminf_{\varepsilon\to
0} \int_{\mathfrak{B}_{\delta,\varepsilon}}\psi^\varepsilon
d\bar{\mathsf{P}},
\end{align*}
where
$
\psi^\varepsilon:=1\wedge\exp\big(-\varepsilon^\kappa\int_0^{T\wedge\tau_{\varepsilon,
\zeta}}\frac{\dot{u}_s-\beta^\varepsilon_s}{(\alpha^\varepsilon_s)^2}d\bar{B}_s
\big).
$
In view of \eqref{A.223}, $ \lim_{\varepsilon\to 0}
\bar{\mathsf{P}}\big(\mathfrak{B}_{\delta,\varepsilon}\big)=1
$
and, also,
$
\lim_{\varepsilon\to 0}\int_\varOmega \psi^\varepsilon d\bar{\mathsf{P}}=1.
$
Hence, owing to $\psi^\varepsilon\le 1$,
\begin{align*}
\int_{\mathfrak{B}_{\delta,\varepsilon}}\psi^\varepsilon d\bar{\mathsf{P}}&=
\int_\varOmega\psi^\varepsilon d\bar{\mathsf{P}}-
\int_{\varOmega\setminus\mathfrak{B}_{\delta,\varepsilon}}\psi^\varepsilon
d\bar{\mathsf{P}}
\\
&\ge
\int_\varOmega\psi^\varepsilon d\bar{\mathsf{P}}-
\int_{\varOmega\setminus\mathfrak{B}_{\delta,\varepsilon}}
d\bar{\mathsf{P}}\xrightarrow[\varepsilon\to 0]{}1,
\end{align*}
that is,
$
\log
\varliminf_{\varepsilon\to
0} \int_{\mathfrak{B}_{\delta,\varepsilon}}\psi^\varepsilon
d\bar{\mathsf{P}}\ge 0.
$

Consequently, by an arbitrariness of $\zeta$, \eqref{LB} is valid.
\end{proof}

\end{document}